\documentclass[12pt]{article}

\usepackage{amssymb}
\usepackage{amsmath}
\usepackage{amscd}
\usepackage[dvips]{graphicx}
\usepackage{mathrsfs}

\numberwithin{equation}{section}

\newtheorem{theorem}{Theorem}[section]

\begin{document}

\title{ON SOJOURN TIMES IN THE $M/M/1$-PS MODEL, CONDITIONED ON THE NUMBER OF OTHER USERS}

\author{
Qiang Zhen\thanks{
Department of Mathematics, Statistics, and Computer Science,
University of Illinois at Chicago, 851 South Morgan (M/C 249),
Chicago, IL 60607-7045, USA.
{\em Email:} qzhen2@uic.edu.}
\and
and
\and
Charles Knessl\thanks{
Department of Mathematics, Statistics, and Computer Science,
University of Illinois at Chicago, 851 South Morgan (M/C 249),
Chicago, IL 60607-7045, USA.
{\em Email:} knessl@uic.edu.\
\newline\indent\indent{\bf Acknowledgement:} This work was partly supported by NSF grant DMS 05-03745 and NSA grant H 98230-08-1-0102.
}}
\date{November 14, 2008 }
\maketitle

\begin{abstract}
\noindent
We consider the $M/M/1$-PS queue with processor sharing. We study the conditional sojourn time distribution of an arriving customer, conditioned on the number of other customers present. A new formula is obtained for the conditional sojourn time distribution, using a discrete Green's function. This is shown to be equivalent to some classic results of Pollaczeck and Vaulot from 1946. Then various asymptotic limits are studied, including large time and/or large number of customers present, and heavy traffic, where the arrival rate is only slightly less than the service rate.
\end{abstract}

\newpage
\section{Introduction}

One of the most interesting service disciplines in queueing theory is that of processor sharing (PS). Here every customer in the system gets an equal fraction of the server or processor, and this has the advantage that shorter jobs get served in less time than, say, under the first-in-first-out (FIFO) discipline.

The PS discipline was introduced by Kleinrock \cite{KLanalysis}, \cite{KLtime}, and has been the subject of much further investigation over the past forty years. In these models one of the main measures of performance is a given (also called tagged) customer's sojourn time distribution, conditioned on the number of other customers in the system upon his arrival. The sojourn time is the total time from when a customer arrives to when that customer leaves the system, after being served.

The $M/M/1$-PS queue assumes Poisson arrivals with rate $\lambda$ and exponential i.i.d. service times with rate $\mu$. The traffic intensity is $\rho=\lambda/\mu$. We shall denote the sojourn time of the tagged customer by $\mathbf{V}$ and the number of other customers present at his arrival instant by $\mathbf{N}$. Then the unconditional sojourn time density is $p(t)dt=\Pr\big[\mathbf{V}\in(t,t+dt)\big]$, while the conditional density, conditioned on $\mathbf{N}$, is $p_n(t)dt=\Pr\big[\mathbf{V}\in(t,t+dt)\big|\mathbf{N}=n]$. For the $M/M/1$-PS model we can remove the conditioning to get
\begin{equation}\label{S1_p_remove}
p(t)=\sum_{n=0}^\infty(1-\rho)\,\rho^np_n(t),
\end{equation}
since $\mathbf{N}$ follows a geometric distribution.

In \cite{COF}, Coffman, Muntz, and Trotter derived an expression for the Laplace transform of the sojourn time distribution, conditioned on both the number seen by an arrival and the amount of service required by the arriving customer, in the $M/M/1$-PS model. Sengupta and Jagerman \cite{SEJA} obtained the moments of the sojourn time distribution conditioned on $\mathbf{N}$, and gave an asymptotic expansion when the number of customers in the system is large. Guillemin and Boyer \cite{GU} formulated $p_n(t)$ as a spectral problem for a self-adjoint operator, and obtained an integral representation for the conditional distribution.

Using the results in \cite{COF}, Morrison \cite{MO} studied the unconditional sojourn time distribution $p(t)$ in the $M/M/1$-PS model, in the heavy traffic limit, where the Poisson arrival rate $\lambda$ is nearly equal to the service rate $\mu$ (thus $\rho=\lambda/\mu\uparrow 1$). Setting $\epsilon=1-\rho$, in \cite{MO} asymptotic results were obtained for the time scales $t=O(1)$, $t=O(\epsilon^{-1})$ and $t=O(\epsilon^{-3})$. Most the mass is concentrated in the range $t=O(\epsilon^{-1})$, and the asymptotic series involves modified Bessel functions.

A service discipline seemingly unrelated to PS is random order service (ROS), where customers are chosen for service at random. The $M/M/1$-ROS model has been studied by many authors, see Vaulot \cite{VA}, Pollaczek \cite{PO}, Riordan \cite{RI}, Kingman \cite{KI} and Flatto \cite{FL}. In \cite{PO} an explicit integral representation is derived for the generating function of the conditional waiting time distribution, from which the following tail behavior of the unconditional waiting time $\mathbf{W}_{\mathrm{ROS}}$ is computed as
\begin{equation}\label{S1_tail}
\Pr\left[\mathbf{W}_{\mathrm{ROS}}>t\right]\sim e^{-\alpha t-\beta t^{1/3}}\gamma t^{-5/6},\textrm{ }\textrm{ }\textrm{ }  t\rightarrow\infty.
\end{equation}
Here $\alpha$, $\beta$ and $\gamma$ are explicitly computed constants, with $\alpha=(1-\sqrt{\rho})^2$. Flatto \cite{FL} obtained an integral representation for the unconditional waiting time distribution and derived the same tail behavior as $t\to\infty$. Cohen \cite{COH} established the following relationship between the sojourn time in the PS model and the waiting time in the ROS model,
\begin{equation}\label{S1_equi}
\Pr\left[\mathbf{V}_{\mathrm{PS}}>t\right]=\frac{1}{\rho}\Pr[\mathbf{W}_{\mathrm{ROS}}>t],
\end{equation}
which extends also to the more general $G/M/1$ case. In  \cite{BO} relations of the form (\ref{S1_equi}) are explored for other models, such as finite capacity queues, repairman problems, and networks.

In this paper we study the conditional sojourn time distribution $p_n(t)$ for the $M/M/1$-PS model in two cases. First we consider a fixed $\rho<1$ and obtain expansions of $p_n(t)$ for $t$ and/or $n\rightarrow\infty$. From these (\ref{S1_tail}) is readily obtained by using (\ref{S1_p_remove}) for $t$ large. Then we consider the heavy traffic limit where $\rho\uparrow 1$, and again obtain approximations for several ranges of the space-time plane. From these results all of the expansions in  \cite{MO} can be recovered. The integral representation in \cite{PO} is used to derive some of the approximations. However, it is difficult to obtain all of the results in this paper from it. Thus, we derive another representation for $p_n(t)$ using a discrete Green's function, which we show to be equivalent to the representation in \cite{PO}.

We mention some related work on various PS models. In \cite{ZH} we studied the sojourn time density conditioned on the service time in the $M/M/1$-PS model for various asymptotic ranges, for both $\rho<1$ and $\rho\approx 1$. The $M/G/1$-PS model was studied by Yashkov \cite{YAproc}, \cite{YAmath}, \cite{YAon} and by Ott \cite{OT}. In \cite{ZW} Zwart and Boxma analyze the $M/G/1$-PS queue with heavy tails, where the service density has algebraic or sub-exponential behavior. Ramaswami \cite{RA} studied the $G/M/1$-PS queue and obtained explicit results for the unconditional moments of the sojourn time. Various asymptotic properties of the conditional and unconditional moments and distribution for this model were derived in \cite{KN}. The $G/G/1$-PS model has not been analyzed exactly, but some approximations are discussed in Sengupta \cite{SEapp} and the tail exponent of the unconditional sojourn time density was derived by Mandjes and Zwart \cite{MA}. A good recent survey of sojourn time asymptotics in PS queues is in Borst, N\'u\~nez-Queija and Zwart \cite{BONU}.

The remainder of the paper is organized as follows. In Section 2 we summarize and briefly discuss our main results (see Theorems 2.1--2.3). In Section 3 we derive the explicit formula for $p_n(t)$ by using a discrete Green's function. In Section 4 we derive the asymptotic results for $p_n(t)$ for moderate traffic intensities $\rho<1$. In Section 5 we consider $p_n(t)$ for $\rho\uparrow 1$, and various scalings of space and time. We discuss a singular perturbation approach to the problem in Section 6.

\newpage

\section{Summary of results}
We consider the $M/M/1$-PS model with arrival rate $\lambda$ and we set the service rate $=1$. Then the traffic intensity is $\rho=\lambda>0$.

It was shown in \cite{COF}, under the stability condition $\rho<1$, that the recurrence equation of the sojourn time density of a tagged customer, conditioned on the number of other customers in the system, is given by
\begin{equation}\label{S2_recu1}
p'_n(t)=\rho\; p_{n+1}(t)-(1+\rho)\; p_n(t)+\frac{n}{n+1}\;p_{n-1}(t),\;\;\;t>0
\end{equation}
with initial condition $p_n(0)=\frac{1}{n+1}$. Taking the Laplace transform of (\ref{S2_recu1}) and multiplying by $n+1$, we have
\begin{equation}\label{S2_recu2}
(n+1)\;\rho\;\widehat{p}_{n+1}(\theta)-(n+1)\;(1+\rho+\theta)\;\widehat{p}_n(\theta)+n\;\widehat{p}_{n-1}(\theta)=-1,
\end{equation}
where $\widehat{p}_n(\theta)=\int^\infty_0 p_n(t)e^{-\theta t}dt$.

Solving the recurrence equation (\ref{S2_recu2}), we obtain the following result.

\begin{theorem} \label{th1}
The Laplace-Stieltjes transform of the conditional sojourn time density has the following form:
\begin{equation}\label{S2_th1_phat}
\widehat{p}_n(\theta)=M\,G_n\sum_{l=0}^n\rho^lH_l+M\,H_n\sum_{l=n+1}^\infty\rho^lG_l,
\end{equation}
where
\begin{equation}\label{S2_th1_M}
M=M(\theta)\equiv z_-\Big(\frac{z_+}{z_-}\Big)^\alpha,
\end{equation}
\begin{equation}\label{S2_th1_G}
G_n=G_n(\theta)\equiv\int_0^{z_-}z^n(z_+-z)^{-\alpha}(z_--z)^{\alpha-1}dz,
\end{equation}
\begin{equation}\label{S2_th1_H}
H_n=H_n(\theta)\equiv\frac{e^{i\alpha\pi}}{2\pi i}\int_{\mathscr{C}}z^n(z_+-z)^{-\alpha}(z-z_-)^{\alpha-1}dz,
\end{equation}
$\mathscr{C}$ is a closed contour in the complex $z$-plane that encircles the segment $[z_-,z_+]$ of the real axis and
\begin{equation}\label{S2_th1_z}
z_\pm=z_\pm(\theta)\equiv\frac{1}{2\rho}\Big[1+\rho+\theta\pm\sqrt{(1+\rho+\theta)^2-4\rho}\Big],
\end{equation}
\begin{equation}\label{S2_th1_alpha}
\alpha=\alpha(\theta)\equiv\frac{z_+}{z_+-z_-}.
\end{equation}

\end{theorem}

We will show in Section 3 that by taking the inverse Laplace transform, the conditional sojourn time density obtained from (\ref{S2_th1_phat}) is equivalent to the result in Pollaczek \cite{PO}:
\begin{eqnarray}\label{S2_polla}
p_n(t) &=& \frac{1}{2\pi i}\oint_{\mathscr{C^\ast}}\frac{1}{z^{n+1}}\bigg[\int_0^\pi e^{-(1+\rho-2\sqrt{\rho}\cos v)t}\frac{(\sqrt{\rho}e^{-iv}-z)^{m_0}}{(\sqrt{\rho}e^{iv}-z)^{m_0+1}}\nonumber\\
       & &\times\; \frac{(1-\sqrt{\rho}e^{iv})^{m_0}}{(1-\sqrt{\rho}e^{-iv})^{m_0+1}}\frac{2\sqrt{\rho}\sin v}{1+e^{\pi \cot v}}dv\bigg]\,dz,
\end{eqnarray}
where
\begin{equation}\label{S2_polla_m0}
m_0=\frac{i}{2}\cot v-\frac{1}{2}
\end{equation}
and the contour $\mathscr{C^\ast}$ is a circle in the complex $z$-plane, centered at the origin and with radius less than $\sqrt{\rho}$.

Using (\ref{S2_th1_phat}) and (\ref{S2_polla}), we obtain the following asymptotic expansions for $p_n(t)$, valid for $\rho<1$ and $n$ and/or $t\to\infty$.

\begin{theorem} \label{th2}
For $\rho<1$, the conditional sojourn time density has the following asymptotic expansions:

\begin{enumerate}

\item $n\to\infty$, $t\to\infty$ with $n/t>1-\rho$,
{\setlength\arraycolsep{2pt}
\begin{eqnarray}\label{S2_th2_1}
p_n(t)&=& \frac{1}{n}\Delta_1^{\frac{\rho}{1-\rho}}+\frac{1}{2(1-\rho)^3n^2}\Delta_1^{\frac{3\rho-2}{1-\rho}}\Big[\rho(2\rho^2+\rho-1)+4\rho^2\Delta_1\log(\Delta_1)\nonumber\\
      &&+\,6\rho(1-\rho)\Delta_1-(\rho^2-\rho+2)\Delta_1^2\Big]+O(n^{-3}),
\end{eqnarray}}
where $\Delta_1\equiv 1-(1-\rho)t/n>0$.

\item $n\to\infty$, $t\to\infty$ with $n/t=1-\rho+O(t^{-1/2})$,
{\setlength\arraycolsep{2pt}
\begin{eqnarray}\label{S2_th2_2}
p_n(t)&\sim&\frac{1}{\sqrt{2\pi}}\;\sqrt{\frac{1-\rho}{1+\rho}}\;n^{-\frac{2-\rho}{2(1-\rho)}}\nonumber\\
&&\times\int_0^\infty y^{\frac{\rho}{1-\rho}}\exp\Big\{-\frac{1-\rho}{2(1+\rho)}\big(y-\Delta_2\big)^2\Big\}dy,
\end{eqnarray}}
where $\Delta_2=\sqrt{n}\,\big[1-(1-\rho)t/n\big]=\sqrt{n}\,\Delta_1=O(1)$.

\item $n\to\infty$, $t\to\infty$ with $0<n/t<1-\rho$,
{\setlength\arraycolsep{2pt}
\begin{eqnarray}\label{S2_th2_3}
p_n(t)&\sim& n^{-1-\frac{1}{2}\sqrt{1+4\rho t^2/n^2}}\;K(\theta_\ast)\exp\bigg\{t\Big(-1-\rho+\sqrt{\frac{n^2}{t^2}+4\rho}\Big)\nonumber\\
&&+\;n\log\Big[\frac{1}{2\rho}\Big(-\frac{n}{t}+\sqrt{\frac{n^2}{t^2}+4\rho}\Big)\Big]\bigg\},
\end{eqnarray}}
where
\begin{equation}\label{S2_th2_3_thetastar}
\theta_\ast=\theta_\ast\Big(\frac{n}{t}\Big)=\sqrt{\frac{n^2}{t^2}+4\rho}-1-\rho,
\end{equation}
\begin{equation}\label{S2_th2_3_K}
K(\theta)=\frac{1}{\sqrt{2\pi}}\,\alpha^\alpha\,\Gamma(\alpha)\,\frac{z_-\,(1-\rho z_-)^{\alpha-1}\big[(1+\rho+\theta)^2-4\rho\big]^{3/4}}{(1-\rho z_+)^{\alpha}\sqrt{1+\rho+\theta}}
\end{equation}
and $\alpha=\alpha(\theta)$ is as in (\ref{S2_th1_alpha}).

\item $n\to\infty$, $t\to\infty$ with $n\:t^{-2/3}\equiv a=O(1)$,
\begin{equation}\label{S2_th2_4}
p_n(t)\sim \rho^{-n/2}\Lambda(n,t)\;\exp\Big\{-(1-\sqrt{\rho})^2t+\Phi(n,t)\Big\}
\end{equation}
where $\Lambda(n,t)$ and $\Phi(n,t)$ have the following expressions in three ranges of $a$ (as shown in Figure 1):

\begin{enumerate}
\item $a\geq (3\sqrt{\rho})^{2/3}$,
\begin{equation}\label{S2_th2_4(1)lambda}
\Lambda(n,t)=\frac{\sqrt{A}\,\exp\Big\{\frac{1+\sqrt{\rho}}{2(1-\sqrt{\rho})}\Big\}}{(1-\sqrt{\rho})\,a^{3/4}\big(3\sqrt{\rho (1+A)}-a^{3/2}\big)^{1/2}\,t},
\end{equation}
\begin{equation}\label{S2_th2_4(1)phi}
\Phi(n,t)=\Big(3\sqrt{\rho}\,\frac{A}{a}-2\sqrt{a(1+A)}\Big)\,t^{1/3},
\end{equation}
where $A=A(a)$ satisfies:
\begin{equation}\label{S2_th2_4(1)A}
2\sqrt{\rho}\,A^{3/2}\,a^{-3/2}=\sqrt{A(1+A)}-\mathrm{arcsinh}\big(\sqrt{A}\big).
\end{equation}

\item $(4\sqrt{\rho}/\pi)^{2/3}< a<(3\sqrt{\rho})^{2/3}$,
\begin{equation}\label{S2_th2_4(2)lambda}
\Lambda(n,t)=\frac{\sqrt{B}\,\exp\Big\{\frac{1+\sqrt{\rho}}{2(1-\sqrt{\rho})}\Big\}}{(1-\sqrt{\rho})\,a^{3/4}\big(-3\sqrt{\rho (1-B)}+a^{3/2}\big)^{1/2}\,t},
\end{equation}
\begin{equation}\label{S2_th2_4(2)phi}
\Phi(n,t)=\Big(-3\sqrt{\rho}\,\frac{B}{a}-2\sqrt{a(1-B)}\Big)\,t^{1/3},
\end{equation}
where $B=B(a)$ satisfies:
\begin{equation}\label{S2_th2_4(2)B}
2\sqrt{\rho}\,B^{3/2}\,a^{-3/2}=-\sqrt{B(1-B)}+\mathrm{arcsin}\big(\sqrt{B}\big).
\end{equation}

\item $0<a\leq(4\sqrt{\rho}/\pi)^{2/3}$,
\begin{equation}\label{S2_th2_4(3)lambda}
\Lambda(n,t)=\frac{\sqrt{C}\,\exp\Big\{\frac{1+\sqrt{\rho}}{2(1-\sqrt{\rho})}\Big\}}{(1-\sqrt{\rho})\,a^{3/4}\big(3\sqrt{\rho (1-C)}+a^{3/2}\big)^{1/2}\,t},
\end{equation}
\begin{equation}\label{S2_th2_4(3)phi}
\Phi(n,t)=\Big(-3\sqrt{\rho}\,\frac{C}{a}+2\sqrt{a(1-C)}\Big)\,t^{1/3},
\end{equation}
where $C=C(a)$ satisfies:
\begin{equation}\label{S2_th2_4(3)C}
2\sqrt{\rho}\,C^{3/2}\,a^{-3/2}=\sqrt{C(1-C)}+\pi-\mathrm{arcsin}\big(\sqrt{C}\big).
\end{equation}

\end{enumerate}

\item $n=O(1)$, $t\to\infty$,
{\setlength\arraycolsep{1pt}
\begin{eqnarray}\label{S2_th2_5}
p_n(t)&\sim& \frac{2^{2/3}3^{-1/2}\pi^{5/6}\rho^{-5/12}}{(1-\sqrt{\rho})\;t^{5/6}} \exp\Big\{-(1-\sqrt{\rho})^2 t-2^{-2/3}3\;\pi^{2/3}\rho^{1/6} t^{1/3}\Big\}\nonumber\\
&\times&\exp\Big(\frac{\sqrt{\rho}}{1-\sqrt{\rho}}\Big)\frac{\rho^{-n/2}}{2\pi i}\oint_{\mathscr{C^\ast}}\frac{1}{(1-z)\;z^{n+1}}\exp\Big(\frac{1}{1-z}\Big)dz.
\end{eqnarray}}

\end{enumerate}
\end{theorem}

We note that for a fixed $a$ and $n,t\to\infty$, $\Phi(n,t)=O(t^{1/3})=O(\sqrt{n})$ and $\Lambda(n,t)=O(t^{-1})=O(n^{-3/2})$. Despite the fact that case 4 has three different expressions, the functions $\Phi$ and $\Lambda$ are smooth along the transition curves $a=n\,t^{-2/3}=(3\sqrt{\rho})^{2/3}$ and $a=(4\sqrt{\rho}/{\pi})^{2/3}$.
We also note that the contour integral in (\ref{S2_th2_5}) is equivalent to the following infinite sum:
\begin{equation*}
\frac{1}{2\pi i}\oint_{\mathscr{C^\ast}}\frac{1}{(1-z)\;z^{n+1}}\exp\Big(\frac{1}{1-z}\Big)dz=\sum_{l=0}^\infty \frac{(n+l)!}{(l!)^2\,n!}.
\end{equation*}

The asymptotic sojourn time density has simpler expressions in some of the matching regions between the scales in Theorem 2.2. We have, between cases 2 and 3,
\begin{eqnarray}\label{S2_th2_matching23}
p_n(t)&\sim&\frac{\sqrt{1-\rho}}{\sqrt{2\pi}}\,\Gamma\Big(\frac{1}{1-\rho}\Big)\Big(1-\rho-\frac{n}{t}\Big)^{-\frac{1}{1-\rho}}\Big(\frac{1+\rho}{n}\Big)^{\frac{1+\rho}{2(1-\rho)}}\nonumber\\
&&\times\exp\Big\{-\frac{t}{2(1+\rho)}\Big(\frac{n}{t}-1+\rho\Big)^2\Big\},
\end{eqnarray}
which is valid for $n\to\infty$, $t\to\infty$, $n/t<1-\rho$ with $O(t^{-1/2})\ll |n/t-1+\rho| \ll 1$.
Between cases 3 and 4(a), we have
\begin{eqnarray}\label{S2_th2_matching34}
p_n(t)&\sim&\frac{\rho^{-n/2}}{\sqrt{2\rho}\,(1-\sqrt{\rho})\,t}\exp\Big\{\frac{1+\sqrt{\rho}}{2(1-\sqrt{\rho})}\Big\}\nonumber\\
&&\times\exp\Big\{-(1-\sqrt{\rho})^2\,t-\frac{n^2}{4\sqrt{\rho}\,t}\Big\}\nonumber\\
&&\times\exp\Big\{\sqrt{\rho}\,\frac{t}{n}\big[\log(\rho)-1+2\log(t)-3\log(n)\big]\Big\},
\end{eqnarray}
which is valid for $t\to\infty$ and $O(t^{2/3})\ll n\ll O(t)$.
Between cases 4(c) and 5, we have
\begin{eqnarray}\label{S2_th2_matching45}
p_n(t)&\sim&\frac{\rho^{-5/12}}{1-\sqrt{\rho}}\,\Big(\frac{\pi}{2}\Big)^{1/3}3^{-1/2}\,\rho^{-n/2}\,t^{-5/6}\,n^{-1/4}\nonumber\\
&&\times\exp\Big\{-(1-\sqrt{\rho})^2\,t-3\Big(\frac{\pi}{2}\Big)^{2/3}\rho^{1/6}\,t^{1/3}\Big\}\nonumber\\
&&\times\exp\Big\{2\sqrt{n}+\frac{1+\sqrt{\rho}}{2(1-\sqrt{\rho})}\Big\},
\end{eqnarray}
which is valid for $t\to\infty$ and $1\ll n\ll O(t^{2/3})$.

By removing the condition on $n$, using (\ref{S2_th2_5}) in (\ref{S1_p_remove}), and noticing the relationship (\ref{S1_equi}) between processor sharing and service in random order, we can recover the results in Pollazcek \cite{PO} and Flatto \cite{FL}, for $p(t)$ as $t\to\infty$.

\vspace{3pt}
We next consider the heavy traffic case, where $\rho$ is close to 1. Letting $\epsilon=1-\rho\to 0^+$, we obtain the following results.

\begin{theorem} \label{th3}
For $\rho=1-\epsilon$ and $\epsilon\to 0^+$, the conditional sojourn time density has the following asymptotic expansions:

\begin{enumerate}

\item $n=O(1)$, $t=O(1)$,
{\setlength\arraycolsep{2pt}
\begin{eqnarray}\label{S2_th3_1}
p_n(t) &\sim& \frac{1}{2\pi i}\oint_{\mathscr{C^\ast}}\frac{1}{z^{n+1}}\bigg[\int_0^\pi e^{-2(1-\cos v)t}\frac{(e^{-iv}-z)^{m_0}}{(e^{iv}-z)^{m_0+1}}\nonumber\\
       & &\times \frac{(1-e^{iv})^{m_0}}{(1-e^{-iv})^{m_0+1}}\frac{2\sin v}{1+e^{\pi \cot v}}dv\bigg]dz,
\end{eqnarray}}
where $m_0$ is defined by (\ref{S2_polla_m0}).

\item $n=O(1)$, $t=\sigma/{\epsilon^3}=O(\epsilon^{-3})$,
{\setlength\arraycolsep{2pt}
\begin{eqnarray}\label{S2_th3_2}
p_n(t) &\sim& 4\sqrt{\pi}\epsilon^{3/2}\frac{u\sqrt{1+4u^2}}{\sqrt{8+3\sigma(1+4u^2)}}\exp\Big[-\big(\frac{1}{8}-\frac{u^2}{2}\big)\sigma-\frac{3+4u^2}{2(1+4u^2)}\Big]\nonumber\\
&&\times\exp\Big\{-\frac{1}{\epsilon}\Big[\big(\frac{1}{4}+3u^2\big)\sigma-\frac{2}{1+4u^2}\Big]\Big\}\nonumber\\
&&\times\frac{1}{2\pi i}\oint_{\mathscr{C^\ast}}\frac{1}{(1-z)\;z^{n+1}}\exp\Big(\frac{1}{1-z}\Big)dz,
\end{eqnarray}}
where $u=u(\sigma)$ satisfies
\begin{equation}\label{S2_th3_2_u}
2u^3\sigma=\pi-\frac{i}{2}\log\Big(\frac{1-2iu}{1+2iu}\Big)+\frac{2u}{1+4u^2}.
\end{equation}

\item $n=\xi/\epsilon=O(\epsilon^{-1})$, $t=\tau/\epsilon=O(\epsilon^{-1})$,
\begin{equation}\label{S2_th3_3}
p_n(t)= \frac{\epsilon}{\xi}e^{-\tau/\xi}+\epsilon^2\Big[\frac{\tau-1}{\xi^2}+\frac{4\tau-\tau^2}{2\xi^3}-\frac{3\tau^2}{2\xi^4}+\frac{\tau^3}{3\xi^5}\Big]e^{-\tau/\xi}+O(\epsilon^3).
\end{equation}

\item $n=\eta/\epsilon^2=O(\epsilon^{-2})$, $t=\sigma/\epsilon^3=O(\epsilon^{-3})$,
\begin{equation}\label{S2_th3_4}
p_n(t)\sim \epsilon^2\,\widetilde{\Lambda}(\eta,\sigma)\exp\Big\{\frac{1}{\epsilon}\Big[ \widetilde{\Phi}(\eta,\sigma)+\frac{\eta}{2}-\frac{\sigma}{4}\Big]\Big\},
\end{equation}
where $\widetilde{\Lambda}(\eta,\sigma)$ and $\widetilde{\Phi}(\eta,\sigma)$ have the following expressions in three ranges of the $(\eta,\sigma)$ plane (as shown in Figure 2):

\begin{enumerate}
\item $\sigma<\frac{1}{3}\eta^{3/2}-\frac{8}{3}$ with $\eta>4$,
{\setlength\arraycolsep{2pt}
\begin{eqnarray}\label{S2_th3_4(1)lambda}
\widetilde{\Lambda}(\eta,\sigma)&=&\frac{2\sqrt{\widetilde{A}\eta(1-4\widetilde{A})}}{\sqrt{(1-4\widetilde{A})^2\Big[3\sigma\eta\sqrt{\eta(1+\widetilde{A}\eta)}-\eta^3\Big]+8\eta\sqrt{\eta(1+\widetilde{A}\eta)}}}\nonumber\\
&&\times\exp\Big[\frac{\eta}{4}-\frac{1}{1-4\widetilde{A}}-\big(\frac{\widetilde{A}}{2}+\frac{1}{8}\big)\sigma\Big],
\end{eqnarray}}
\begin{equation}\label{S2_th3_4(1)phi}
\widetilde{\Phi}(\eta,\sigma)=3\widetilde{A}\sigma-2\sqrt{\eta(1+\widetilde{A}\eta)}+\frac{2}{1-4\widetilde{A}},
\end{equation}
where $\widetilde{A}=\widetilde{A}(\eta,\sigma)$ satisfies
{\setlength\arraycolsep{2pt}
\begin{eqnarray}\label{S2_th3_4(1)A}
2\widetilde{A}^{3/2}\sigma&=& -\frac{2\sqrt{\widetilde{A}}}{1-4\widetilde{A}}-\mathrm{arcsinh}\Big(\sqrt{\widetilde{A}\eta}\Big)+\mathrm{arcsinh}\bigg(\sqrt{\frac{4\widetilde{A}}{1-4\widetilde{A}}}\bigg)\nonumber\\
&&+\;\sqrt{\widetilde{A}\eta(1+\widetilde{A}\eta)}.
\end{eqnarray}}

\item $\frac{1}{3}\eta^{3/2}-\frac{8}{3}\leq\sigma\leq\frac{1}{2}\eta^{3/2}\Big[\frac{\pi}{2}+\frac{4\sqrt{\eta}}{4+\eta}-\arcsin\Big(\sqrt{\frac{4}{4+\eta}}\Big)\Big],\;\sigma>0$,
{\setlength\arraycolsep{2pt}
\begin{eqnarray}\label{S2_th3_4(2)lambda}
\widetilde{\Lambda}(\eta,\sigma)&=&\frac{2\sqrt{\widetilde{B}\eta(1+4\widetilde{B})}}{\sqrt{(1+4\widetilde{B})^2\Big[\eta^3-3\sigma\eta\sqrt{\eta(1-\widetilde{B}\eta)}\Big]-8\eta\sqrt{\eta(1-\widetilde{B}\eta)}}}\nonumber\\
&&\times\exp\Big[\frac{\eta}{4}-\frac{1}{1+4\widetilde{B}}-\big(-\frac{\widetilde{B}}{2}+\frac{1}{8}\big)\sigma\Big],
\end{eqnarray}}
\begin{equation}\label{S2_th3_4(2)phi}
\widetilde{\Phi}(\eta,\sigma)=-3\widetilde{B}\sigma-2\sqrt{\eta(1-\widetilde{B}\eta)}+\frac{2}{1+4\widetilde{B}},
\end{equation}
where $\widetilde{B}=\widetilde{B}(\eta,\sigma)$ satisfies
{\setlength\arraycolsep{2pt}
\begin{eqnarray}\label{S2_th3_4(2)B}
2\widetilde{B}^{3/2}\sigma&=& \frac{2\sqrt{\widetilde{B}}}{1+4\widetilde{B}}+\arcsin\Big(\sqrt{\widetilde{B}\eta}\Big)-\arcsin\bigg(\sqrt{\frac{4\widetilde{B}}{1+4\widetilde{B}}}\bigg)\nonumber\\
&&-\;\sqrt{\widetilde{B}\eta(1-\widetilde{B}\eta)}.
\end{eqnarray}}

\item $\sigma> \frac{1}{2}\eta^{3/2}\Big[\frac{\pi}{2}+\frac{4\sqrt{\eta}}{4+\eta}-\arcsin\Big(\sqrt{\frac{4}{4+\eta}}\Big)\Big]$,
{\setlength\arraycolsep{2pt}
\begin{eqnarray}\label{S2_th3_4(3)lambda}
\widetilde{\Lambda}(\eta,\sigma)&=&\frac{2\sqrt{\widetilde{C}\eta(1+4\widetilde{C})}}{\sqrt{(1+4\widetilde{C})^2\Big[3\sigma\eta\sqrt{\eta(1-\widetilde{C}\eta)}+\eta^3\Big]+8\eta\sqrt{\eta(1-\widetilde{C}\eta)}}}\nonumber\\
&&\times\exp\Big[\frac{\eta}{4}-\frac{1}{1+4\widetilde{C}}-\big(-\frac{\widetilde{C}}{2}+\frac{1}{8}\big)\sigma\Big],
\end{eqnarray}}
\begin{equation}\label{S2_th3_4(3)phi}
\widetilde{\Phi}(\eta,\sigma)=-3\widetilde{C}\sigma+2\sqrt{\eta(1-\widetilde{C}\eta)}+\frac{2}{1+4\widetilde{C}},
\end{equation}
where $\widetilde{C}=\widetilde{C}(\eta,\sigma)$ satisfies
{\setlength\arraycolsep{2pt}
\begin{eqnarray}\label{S2_th3_4(3)C}
2\widetilde{C}^{3/2}\sigma&=& \pi+\frac{2\sqrt{\widetilde{C}}}{1+4\widetilde{C}}-\arcsin\Big(\sqrt{\widetilde{C}\eta}\Big)-\arcsin\bigg(\sqrt{\frac{4\widetilde{C}}{1+4\widetilde{C}}}\bigg)\nonumber\\
&&+\;\sqrt{\widetilde{C}\eta(1-\widetilde{C}\eta)}.
\end{eqnarray}}

\end{enumerate}
\end{enumerate}
\end{theorem}

In the heavy traffic case we can also get much more explicit expressions in the matching regions. Between the regions 1 and 2 in Theorem 2.3,
\begin{eqnarray}\label{S2_th3_matching12}
p_n(t)&\sim& \frac{2\sqrt{\pi}}{\sqrt{3\,t}}\exp\Big\{-2^{-4/3}\,3\,\pi^{2/3}\,t^{1/3}-\frac{1}{2}\Big\}\nonumber\\
&&\times\frac{1}{2\pi i}\oint_{\mathscr{C^\ast}}\frac{1}{(1-z)\;z^{n+1}}\exp\Big({\frac{1}{1-z}}\Big)dz,
\end{eqnarray}
which is valid for $n=O(1)$ and $1\ll t\ll O(\epsilon^{-2})$.
In the matching region between cases 2 and 4(c), where $1\ll n\ll O(\epsilon^{-2})$ and $t=O(\epsilon^{-3})$,
\begin{eqnarray}\label{S2_th3_matching24}
p_n(t) &\sim& \frac{2\,\epsilon^2\,\eta^{-1/4}\,u\,\sqrt{1+4u^2}}{\sqrt{8+3\sigma(1+4u^2)}}\exp\Big\{-\big(\frac{1}{8}-\frac{u^2}{2}\big)\sigma-\frac{3+4u^2}{2(1+4u^2)}+\frac{1}{2}\Big\}\nonumber\\
&&\times\;\exp\Big\{-\frac{1}{\epsilon}\Big[\big(\frac{1}{4}+3u^2\big)\sigma-\frac{2}{1+4u^2}-2\sqrt{\eta}\Big]\Big\},
\end{eqnarray}
and $u=u(\sigma)$ satisfies (\ref{S2_th3_2_u}).

By removing the condition on $n$, using the results (\ref{S2_th3_3}) and (\ref{S2_th3_4}) in (\ref{S1_p_remove}), we can recover the results for $p(t)$ in Morrison \cite{MO} for the time ranges $t=O(\epsilon^{-1})$ and $t=O(\epsilon^{-3})$.

Using our results for $p_n(t)$ we can also obtain some conditional limit laws for $p(n|t)=p_n(t)\,(1-\rho)\,\rho^n/p(t)$, which is the conditional probability of finding $n$ other customers in the system, given the tagged customer's sojourn time. For $t=\tau/\epsilon=O(\epsilon^{-1})$ most of the mass occurs in the range $n=\xi/\epsilon=O(\epsilon^{-1})$ and we have
\begin{equation}\label{S2_limit1}
p(n|t)\approx \frac{\epsilon}{2\xi K_0(2\sqrt{\tau})}e^{-\xi-\tau/\xi},
\end{equation}
where $K_0(\cdot)$ is the modified Bessel function. For $\tau\to\infty$ this simplifies to the Gaussian
\begin{equation}\label{S2_limit1_Gaussian}
p(n|t)\approx \frac{\epsilon}{\sqrt{\pi}\,\tau^{1/4}}\exp\Big\{-\frac{1}{\sqrt{\tau}}(\xi-\sqrt{\tau})^2\Big\},
\end{equation}
which applies for $\xi=\sqrt{\tau}+O(\tau^{1/4})$.

When $n=\eta/\epsilon^2$ and $t=\sigma/\epsilon^3$, we obtain from case 4(c) in Theorem 2.3
\begin{equation}\label{S2_limit2}
p(n|t)\approx \epsilon^{3/2}\sqrt{\frac{Q(\sigma)}{2\pi}}\exp\bigg\{-\frac{Q(\sigma)}{2\epsilon}\Big(\eta-\frac{4}{4\widehat{C}+1}\Big)^2\bigg\},
\end{equation}
where $\widehat{C}=\widehat{C}(\sigma)$ satisfies
\begin{equation}\label{S2_limit2_C}
2\widehat{C}^{3/2}\sigma=\pi+\frac{4\sqrt{\widehat{C}}}{1+4\widehat{C}}-2\arcsin\bigg(\sqrt{\frac{4\widehat{C}}{1+4\widehat{C}}}\bigg)
\end{equation}
and
\begin{equation}\label{S2_limit2_Q}
Q(\sigma)=\frac{(1+4\widehat{C})^2[3\sigma(1+4\widehat{C})^2+16]}{16[3\sigma(1+4\widehat{C})^2+16(1+2\widehat{C})]}.
\end{equation}
The Gaussian limit law in (\ref{S2_limit2}) applies for $\eta=4/[4\widehat{C}(\sigma)+1]+O(\sqrt{\epsilon})$.

\newpage

\section{Brief derivation of Theorem 2.1}

We use a discrete Green's function to derive (\ref{S2_th1_phat}). Consider the recurrence equation (\ref{S2_recu2}). The discrete Green's function $\mathscr{G}(\theta;n,l)$ satisfies
\begin{eqnarray}\label{S3_dgf}
&&(n+1)\rho\,\mathscr{G}(\theta;n+1,l)-(n+1)(1+\rho+\theta)\,\mathscr{G}(\theta;n,l)\nonumber\\
&&\quad\quad\quad\quad\quad\quad+n\,\mathscr{G}(\theta;n-1,l)=-\delta(n,l),\quad (n,l\geq 0)
\end{eqnarray}
where $\delta(n,l)=1_{\{n=l\}}$ is the Kronecker delta.
To construct the Green's function requires that we have two linearly independent solutions to
\begin{equation}\label{S3_dgf_Homo}
(n+1)\rho\,G(\theta;n+1,l)-(n+1)(1+\rho+\theta)\,G(\theta;n,l)+n\,G(\theta;n-1,l)=0,
\end{equation}
which is the homogeneous version of (\ref{S3_dgf}).

We seek solutions of (\ref{S3_dgf_Homo}) in the form
$$ G_n=\int_\mathscr{D}z^ng(z)dz,$$
where the function $g(z)$ and path $\mathscr{D}$ of integration in the complex $z$-plane are to be determined. Using in the above in (\ref{S3_dgf_Homo}) and integrating by part yields
\begin{eqnarray}\label{S3_IBP}
&&z^ng(z)\big[\rho z^2-(1+\rho+\theta)z+1\big]\Big|_\mathscr{D}\nonumber\\
&&\quad\quad -\int_\mathscr{D}z^n\big[(\rho z^2-(1+\rho+\theta)z+1)g'(z)+\rho zg(z)\big]dz=0.
\end{eqnarray}
The first term represents contributions from the endpoints of the contour $\mathscr{D}$.

If (\ref{S3_IBP}) is to hold for all $n$ the integrand must vanish, so that $g(z)$ must satisfy the differential equation
\begin{equation}\label{S3_gz}
\big[\rho z^2-(1+\rho+\theta)z+1\big]g'(z)+\rho zg(z)=0.
\end{equation}
We denote the roots of $\rho z^2-(1+\rho+\theta)z+1=0$ by $z_+$ and $z_-$, (with $z_+>z_->0$ for real $\theta$). These are given by (\ref{S2_th1_z}) and if $\alpha$ is defined by (\ref{S2_th1_alpha}), the solution for $g(z)$ is
$$ g(z)=(z_+-z)^{-\alpha}(z_--z)^{\alpha-1}.$$

If the path of integration $\mathscr{D}$ is chosen as the segment $[0,z_-]$ of the real axis, then (\ref{S3_IBP}) is satisfied for $n\geq 1$. Thus, we have $G_n$ as in (\ref{S2_th1_G}). We note that $G_n$ decays as $n\to\infty$, and is asymptotically given by
\begin{equation}\label{S3_Gasymp}
G_n\sim\frac{\Gamma(\alpha)}{n^\alpha}z_-^{\alpha+n}(z_+-z_-)^{-\alpha},\quad n\to\infty.
\end{equation}
However, $G_n$ becomes infinite as $n\to-1$, which means that $nG_{n-1}$ goes to a nonzero limit as $n\to 0$. Thus $G_n$ is not an acceptable solution to (\ref{S3_dgf_Homo}) at $n=0$.

To construct a second solution to (\ref{S3_dgf_Homo}), we consider another path of integration, $\mathscr{C}$, which is a closed contour in the complex $z$-plane, around the segment $(z_-,z_+)$ of the real axis. Then (\ref{S3_IBP}) is again satisfied as the endpoint contributions from both $z=z_-$ and $z=z_+$ vanish. Thus, we have another solution of (\ref{S3_dgf_Homo}), $H_n$, which is given by (\ref{S2_th1_H}). $H_n$ is finite as $n\to-1$, but grows as $n\to\infty$:
\begin{equation}\label{S3_Hasymp}
H_n\sim\frac{n^{\alpha-1}}{\Gamma(\alpha)}z_+^{n+1-\alpha}(z_+-z_-)^{\alpha-1},\quad n\to\infty.
\end{equation}

Thus, the discrete Green's function can be represented by
\begin{eqnarray}\label{S3_G_step}
\mathscr{G}(\theta;n,l) = \left\{ \begin{array}{ll}
H_l\,G_n\,\mathscr{G}_0 & \textrm{if $n\geq l$}\\
G_l\,H_n\,\mathscr{G}_0 & \textrm{if $0\leq n<l$},
\end{array} \right.
\end{eqnarray}
which has acceptable behavior both at $n=0$ and as $n\to\infty$. Here $\mathscr{G}_0$ depends only upon $\theta$ and $l$.

To determine $\mathscr{G}_0$, we let $n=l$ in (\ref{S3_dgf_Homo}) and use the identities
$$(l+1)\,\rho\,H_{l+1}=(l+1)\,(1+\rho+\theta)\,H_l-l\,H_{l-1},$$
$$(l+1)\,\rho\,G_{l+1}=(l+1)\,(1+\rho+\theta)\,G_l-l\,G_{l-1}.$$
From the above we can infer a simple difference equation for the discrete Wronskian $G_l\,H_{l+1}-G_{l+1}\,H_l$, whose solution we write as
\begin{equation}\label{S3_G1}
G_l\,H_{l+1}-G_{l+1}\,H_l=\frac{1}{(l+1)\,\rho\,^l\,\mathscr{G}_1},
\end{equation}
where $\mathscr{G}_1=\mathscr{G}_1(\theta)$ depends upon $\theta$ only. Then using (\ref{S3_G_step}) in (\ref{S3_dgf_Homo}) with $n=l$ shows that $\mathscr{G}_0$ and $\mathscr{G}_1$ are related by $\mathscr{G}_0=\rho^{l-1}\,\mathscr{G}_1$.

Letting $l\to\infty$ in (\ref{S3_G1}) and using the asymptotic results (\ref{S3_Gasymp}) and (\ref{S3_Hasymp}), we determine $\mathscr{G}_1$ and then obtain
$$\mathscr{G}_0=\rho^lz_-\Big(\frac{z_+}{z_-}\Big)^\alpha.$$

Then, we multiply (\ref{S3_dgf}) by the solution $\widehat{p}_l(\theta)$ to (\ref{S2_recu2}) and sum over all $l\geq 0$. After some manipulation this yield
$$\widehat{p}_n(\theta)=\sum_{l=0}^\infty\mathscr{G}(\theta;n,l),$$
which is equivalent to (\ref{S2_th1_phat}).

The inverse Laplace transform gives the conditional sojourn time density $p_n(t)$ as
\begin{equation}\label{S3_p}
p_n(t)=\frac{1}{2\pi i}\int_{Br}\widehat{p}_n(\theta)e^{\theta t}d\theta,
\end{equation}
where $Br$ is a vertical contour in the complex $\theta$-plane, with $\Re(\theta)>-(1-\sqrt{\rho})^2$.

Now we show the equivalence between (\ref{S3_p}) and (\ref{S2_polla}). We rewrite (\ref{S2_th1_phat}) as
\begin{equation}\label{S3_phat_rewrite}
\widehat{p}_n(\theta)=M\,H_n\sum_{l=0}^\infty\rho^l\,G_l+M\,\sum_{l=0}^n\rho^l\big[G_n\,H_l-G_l\,H_n\big].
\end{equation}
We deform the contour of integration in (\ref{S3_p}) and evaluate the integrand along the line segments just above and just below the branch cut $\Re(\theta)\in[-(1+\sqrt{\rho})^2, -(1-\sqrt{\rho})^2]$. We denote these values of $\widehat{p}_n(\theta)$ by $\Psi_1(\theta)$ and $\Psi_2(\theta)$, respectively. Then (\ref{S3_p}) becomes
\begin{equation}\label{S3_p_psi}
p_n(t)=\frac{1}{2\pi i}\int_{-(1+\sqrt{\rho})^2}^{-(1-\sqrt{\rho})^2}\big[\Psi_2(\theta)-\Psi_1(\theta)\big]e^{\theta t}d\theta
\end{equation}
and we note that $\Psi_1(\theta)$ changes to $\Psi_2(\theta)$ after making the transformation $z_+\to z_-$ and $\alpha\to 1-\alpha$.

We evaluate $H_n$ in (\ref{S2_th1_H}) by branch cut integration, which yields, for $0<\alpha<1$,
\begin{eqnarray}\label{S3_equi_Hn}
H_n&=&\frac{\sin{\alpha\pi}}{\pi}\int_{z_-}^{z_+}\xi^n\,(\xi-z_-)^{\alpha-1}(z_+-\xi)^{-\alpha}d\xi\\
&=&z_-^n\;_2F_1\big(\alpha,-n;1;\frac{1}{1-\alpha}\big),\nonumber
\end{eqnarray}
where $_2F_1$ is the hypergeometric function. Then we observe that $H_n$ and $M$ are both invariant under the transformation $z_+\to z_-$ and $\alpha\to 1-\alpha$. Thus we have
\begin{eqnarray}\label{S3_Psi1}
\Psi_1(\theta)&=&M\,H_n\int_0^{z_-}\frac{1}{1-\rho z}(z_+-z)^{-\alpha}(z_--z)^{\alpha-1}dz\nonumber\\
&+&M\sum_{l=0}^n\rho^l\int_0^{z_-}(z_+-z)^{-\alpha}(z_--z)^{\alpha-1}\big[z^nH_l-z^lH_n\big]dz.
\end{eqnarray}
Here we used (\ref{S2_th1_G}) to evaluate (\ref{S3_phat_rewrite}) above the branch cut. But, $\Psi_2(\theta)$ is the same as $\Psi_1(\theta)$, except for changing the upper limits on both of the integrals in (\ref{S3_Psi1}) from $z_-$ to $z_+$. Also, the function $(z_+-z)^{-\alpha}(z_--z)^{\alpha-1}$ is invariant under the map $z_+\to z_-$, $\alpha\to 1-\alpha$. The difference is $\Psi_2-\Psi_1$, thus
{\setlength\arraycolsep{2pt}
\begin{eqnarray}\label{S3_psi_diff}
\Psi_2(\theta)&-&\Psi_1(\theta)=M\,H_n\int_{z_-}^{z_+}\frac{1}{1-\rho z}(z_+-z)^{-\alpha}(z_--z)^{\alpha-1}dz\nonumber\\
&&+\,M\Big[\sum_{l=0}^n\rho^l\,H_l\int_{z_-}^{z_+}(z_+-z)^{-\alpha}(z_--z)^{\alpha-1}dz\nonumber\\
&&-\,H_n\int_{z_-}^{z_+}\frac{1-(\rho z)^{n+1}}{1-\rho z}(z_+-z)^{-\alpha}(z_--z)^{\alpha-1}dz\Big].
\end{eqnarray}}

Using (\ref{S3_equi_Hn}) and after some calculation, we find that the second part in (\ref{S3_psi_diff}) is zero, and the integral in the first part can be evaluated by using contour integration (using the fact that there is a simple pole at $z=1/\rho$). Thus, (\ref{S3_p_psi}) becomes
\begin{eqnarray}\label{S3_p_trans}
p_n(t)&=&\frac{1}{2\pi i}\int_{\mathscr{C}}\int_{-(1+\sqrt{\rho})^2}^{-(1-\sqrt{\rho})^2}e^{\theta\,t}M\,z^n\,(z_+-z)^{-\alpha}(z-z_-)^{\alpha-1}\nonumber\\
&&\times \frac{e^{i\alpha\pi}}{1-e^{2\pi i\alpha}}\Big(\frac{z_-}{1-z_-}\Big)^\alpha\Big(\frac{1-z_+}{z_+}\Big)^{\alpha-1} d\theta\, dz.
\end{eqnarray}
Using the transformations $\theta\to-1-\rho+2\sqrt{\rho}\cos v$ in (\ref{S3_p_trans}) and $z\to 1/z$ in (\ref{S2_th1_H}), and changing the order of integration, we see the equivalence between (\ref{S3_p}) and (\ref{S2_polla}). Note that with this transformation, $(\alpha, z_+, z_-)$ becomes $(-m_0, e^{iv}/\sqrt{\rho}, e^{-iv}/\sqrt{\rho})$.

\newpage

\section{Asymptotic results for the case $\rho<1$}

We assume that the traffic intensity $\rho$ is fixed and less than one. We sketch the main points in deriving Theorem 2.2. We first consider $n,t\to\infty$ with $n/t>1-\rho$ and use the result in (\ref{S2_th1_phat}). From (\ref{S3_Gasymp}) and (\ref{S3_Hasymp}), we notice that the first term in (\ref{S2_th1_phat}) dominates the second, and thus the Laplace transform is asymptotically given by
\begin{eqnarray}
\widehat{p}_n(\theta)&\sim& M\,G_n\,\sum_{l=0}^{n}\rho^l\,H_l\nonumber\\
&\sim&\frac{1}{\rho\,(z_+-z_-)}\sum_{l=0}^n\frac{l^{\alpha-1}}{n^\alpha}z_-^{n-l}\label{S4_th2_1_sum}\\
&\sim&\frac{1}{\rho\,(z_+-z_-)}\int_0^1y^{\alpha-1}z_-^{n(1-y)}dy.\label{S4_th2_1_int}
\end{eqnarray}
Here we used the Euler-Maclaurin sum formula to approximate the sum in (\ref{S4_th2_1_sum}) by an integral. By scaling $\theta=\nu/n=O(1/n)$ and using 
$$z_-=1-\frac{\theta}{1-\rho}+O(\theta^2)\quad \textrm{as}\quad \theta\to 0,$$
(\ref{S4_th2_1_int}) becomes
\begin{eqnarray}\label{S4_th2_1_truncate}
\widehat{p}_n(\theta)&\sim& \frac{1}{1-\rho}\int_0^{1}y^{\frac{\rho}{1-\rho}}\,e^{-\frac{1-y}{1-\rho}\nu}dy\nonumber\\
&=&\int_0^{\frac{1}{1-\rho}}\big[1-(1-\rho)\,y\big]^{\frac{\rho}{1-\rho}}e^{-\nu\,y}\,dy.
\end{eqnarray}
Then we multiply $\widehat{p}_n(\theta)$ by $e^{\theta\,t}d\theta=n^{-1}\,e^{\nu\,t/n}\,d\nu$ and invert the transform to obtain
\begin{equation}\label{S4_th2_1_pleading}
p_n(t)\sim\frac{1}{n}\Big[1-(1-\rho)\frac{t}{n}\Big]^{\frac{\rho}{1-\rho}}, \quad \frac{n}{t}>1-\rho.
\end{equation}
To obtain the second term in (\ref{S2_th2_1}), we need the correction terms in (\ref{S4_th2_1_int}), for which we also need the second terms in the approximations in (\ref{S3_Gasymp}) and (\ref{S3_Hasymp}). 

This analysis suggests that $p_n(t)$ is approximately zero in the range $0<n/t<1-\rho$. We shall show that in this sector the density is exponentially small. But, we first investigate the case $(1-\rho)t\approx n$.

Thus, we consider $n,t\to\infty$ with $n/t=1-\rho+O(t^{-1/2})$. We can still use (\ref{S4_th2_1_sum}) but now scale $l=y\,\sqrt{n}=O(\sqrt{n})$, and approximate the sum in (\ref{S4_th2_1_sum}) by 
\begin{equation*}
\widehat{p}_n(\theta)\sim\frac{z_-^n}{\rho\,(z_+-z_-)\,n^{\alpha/2}}\int_0^\infty y^{\alpha-1}\,e^{-\sqrt{n}\log(z_-)\,y}dy.
\end{equation*}
Taking the inverse Laplace transform and scaling $\theta=\varpi/\sqrt{n}=O(1/\sqrt{n})$, we have
\begin{eqnarray}\label{S4_th2_2}
p_n(t)&\sim&\frac{1}{2\pi i}\int_{Br}\frac{z_-^n}{\rho(z_+-z_-)n^{\alpha/2}}\Big(\int_0^\infty y^{\alpha-1}e^{-\sqrt{n}\log(z_-)\,y}dy\Big) e^{\theta t}d\theta\nonumber\\
&\sim& \frac{1}{1-\rho}\,n^{-\frac{2-\rho}{2(1-\rho)}}\int_0^\infty y^\frac{\rho}{1-\rho}\Big(\frac{1}{2\pi i}\int_{Br}e^{F(\varpi,\,y)}d\varpi\Big)dy,
\end{eqnarray}
where
\begin{equation*}
F(\varpi,y)=\frac{1+\rho}{2(1-\rho)^3}\,\varpi^2+\Big(\frac{t}{\sqrt{n}}-\frac{\sqrt{n}}{1-\rho}+\frac{y}{1-\rho}\Big)\varpi.
\end{equation*}
Then by using the identity
\begin{equation*}
\frac{1}{2\pi i}\int_{Br}e^{C_0\varpi^2+C_1\varpi}d\varpi=\frac{1}{2\sqrt{\pi\,C_0}}\exp\Big(-\frac{C_1^2}{4C_0}\Big),
\end{equation*}
and noting that 
$$\frac{t}{\sqrt{n}}-\frac{\sqrt{n}}{1-\rho}=-\frac{\Delta_2}{1-\rho},$$
we explicitly evaluate the integral over $\varpi$ in (\ref{S4_th2_2}) to obtain (\ref{S2_th2_2}).

For the case $n,t\to\infty$ with $0<n/t<1-\rho$, we rewrite (\ref{S2_th1_phat}) as
\begin{eqnarray}\label{S4_th2_3}
\widehat{p}_n(\theta)&=&M\,G_n\sum_{l=0}^\infty\rho^l\,H_l\;+\;M\sum_{l=n+1}^\infty\rho^l\,(H_n\,G_l-H_l\,G_n)\nonumber\\
&\sim& M\,G_n\sum_{l=0}^\infty\rho^l\,H_l.
\end{eqnarray}
The second sum is negligible in view of (\ref{S3_Gasymp}) and (\ref{S3_Hasymp}), and the fact that we will have $\theta<0$ on this scale. The sum in (\ref{S4_th2_3}) can be calculated exactly by using (\ref{S2_th1_H}), contour integration and the residue theorem, which yields
\begin{equation}\label{S4_th3_Hsum}
\sum_{l=0}^\infty\rho^l\,H_l=\frac{1}{\rho}\,\Big(\frac{1}{\rho}-z_+\Big)^{-\alpha}\Big(\frac{1}{\rho}-z_-\Big)^{\alpha-1}=\frac{(1-\rho z_-)^{\alpha-1}}{(1-\rho z_+)^\alpha}.
\end{equation}
Here we used the identity $\rho\,z_+\,z_-=1$. Using (\ref{S3_Gasymp}) and (\ref{S4_th3_Hsum}) in (\ref{S4_th2_3}), then taking the inverse Laplace transform, we have
\begin{equation*}
p_n(t)\sim\frac{1}{2\pi i}\int_{Br}\frac{\Gamma(\alpha)}{n^{\alpha}}\frac{z_-\,z_+^\alpha\,(1-\rho\,z_-)^{\alpha-1}}{(z_+-z_-)^\alpha\,(1-\rho\,z_+)^\alpha}\,e^{\widetilde{F}(\theta)}d\theta,
\end{equation*}
where $\widetilde{F}(\theta)=\theta t+n\log(z_-(\theta))$. There is a saddle point at $\theta=\theta_\ast<0$ which satisfies $\widetilde{F}'(\theta)=t+n\,z_-'(\theta)/z_-(\theta)=0$, and this leads to $\theta_\ast$ in (\ref{S2_th2_3_thetastar}). Hence using the saddle point method leads to (\ref{S2_th2_3}).

The expression (\ref{S2_th2_matching23}) in the matching region between cases 2 and 3 follows by letting $\Delta_2\to -\infty$ in (\ref{S2_th2_2}), or letting $n/t\to 1-\rho$ in (\ref{S2_th2_3}) (which corresponds to $\theta_\ast\to 0$).

Now we consider $n,t\to\infty$ with $n\,t^{-2/3}\equiv a=O(1)$. We first note that $\theta_\ast\to\theta_c=-(1-\sqrt{\rho})^2$ as $n/t\to 0$ in (\ref{S2_th2_3_thetastar}). Expressions (\ref{S4_th2_3}) and (\ref{S4_th3_Hsum}) are still valid and we have
\begin{equation*}
\widehat{p}_n(\theta)\sim\frac{z_-}{1-\rho\,z_-}\Big(\frac{z_+-1}{z_--1}\Big)^\alpha\,G_n.
\end{equation*}
Then, by taking the inverse Laplace transform, the conditional sojourn time density is asymptotically given by the double integral
\begin{equation}\label{S4_th2_4(1)pnt}
p_n(t)\sim\frac{1}{2\pi i}\int_{Br}\frac{z_-}{1-\rho\,z_-}\Big(\frac{z_+-1}{z_--1}\Big)^\alpha\,\Big(\int_0^{z_-}z^n\,\frac{(z_--z)^{\alpha-1}}{(z_+-z)^{\alpha}}\,dz\Big)\,e^{\theta\,t}\,d\theta.
\end{equation}
Scaling $\theta=\theta_c+s/n$ (with $s\geq 0$), we notice that 
\begin{equation*}
z_\pm=\rho^{-1/2}\pm\rho^{-3/4}\,\frac{\sqrt{s}}{\sqrt{n}}+O\Big(\frac{1}{n}\Big).
\end{equation*}
Thus, we scale $z=z_--\rho^{-3/4}\,y/\sqrt{n}$ (with $y>0$). The inner integral in (\ref{S4_th2_4(1)pnt}), which is $G_n$, is asymptotically equal to 
\begin{equation}\label{S4_th2_4(1)Gn}
G_n\sim\int_0^\infty\frac{\rho^{-n/2}}{\sqrt{y(y+2\sqrt{s})}}\exp\Big\{\sqrt{n}\,\phi(y,s)+\frac{s}{2\sqrt{\rho}}-\frac{(y+\sqrt{s})^2}{2\sqrt{\rho}}\Big\}\,dy,
\end{equation}
where
\begin{equation*}
\phi(y,s)=\frac{\rho^{1/4}}{2\sqrt{s}}\log\Big(\frac{y}{y+2\sqrt{s}}\Big)-\rho^{-1/4}\,(y+\sqrt{s}).
\end{equation*}
The function $\phi(y,s)$ has its maximum at $y=y_\ast(s)=\sqrt{s+\sqrt{\rho}}-\sqrt{s}>0$, which satisfies $\phi_y(y_\ast,s)=0$ and $\phi_{yy}(y_\ast,s)<0$. Hence, using the Laplace method in (\ref{S4_th2_4(1)Gn}), we have 
\begin{equation}\label{S6_S4_th2_4(1)Gnasym}
G_n\sim \frac{\sqrt{\pi}}{\big[n\,(s+\sqrt{\rho})\big]^{1/4}}\,\rho^{-n/2}\,\rho^{1/8}\,e^{-1/2}\,e^{\sqrt{n}\,\phi(y_\ast,s)}.
\end{equation} 
In the limit of $\theta\to \theta_c$ we have
\begin{equation*}
\frac{z_-}{1-\rho\,z_-}\Big(\frac{z_+-1}{z_--1}\Big)^\alpha\sim\frac{1}{\sqrt{\rho}(1-\sqrt{\rho})}\exp\Big(\frac{1}{1-\sqrt{\rho}}\Big),
\end{equation*}
with which $p_n(t)$ in (\ref{S4_th2_4(1)pnt}) asymptotically becomes
\begin{eqnarray}\label{S4_th2_4(1)pnt2}
p_n(t)&\sim&\frac{\sqrt{\pi}\,\rho^{-n/2-3/8}}{(1-\sqrt{\rho})\,n^{5/4}}\exp\Big\{\frac{1+\sqrt{\rho}}{2(1-\sqrt{\rho})}\Big\}\,e^{-(1-\sqrt{\rho})^2\,t}\nonumber\\
&&\times\,\frac{1}{2\pi i}\int_{Br'}(s+\sqrt{\rho})^{-1/4}\,e^{\sqrt{n}\,\Phi_0(s)}ds,
\end{eqnarray}
where $\Phi_0(s)=\phi(y_\ast(s),s)+a^{-3/2}\,s$, and $Br'$ is a vertical contour in the complex $s$-plane. The saddle point $s_\ast=s_\ast(a)$ satisfies $\Phi_0'(s_\ast)=0$, which implies that
\begin{equation}\label{S4_th2_4_sstar}
4\,\rho^{1/4}\,a^{-3/2}\,s_\ast^{3/2}-\sqrt{\rho}\,\log\bigg(\frac{\sqrt{s_\ast+\sqrt{\rho}}-\sqrt{s_\ast}}{\sqrt{s_\ast+\sqrt{\rho}}+\sqrt{s_\ast}}\bigg)-2\sqrt{s_\ast\big(s_\ast+\sqrt{\rho}\big)}=0.
\end{equation}
If let $s_\ast(a)=\sqrt{\rho}\,A\geq 0$, (\ref{S4_th2_4_sstar}) is equivalent to (\ref{S2_th2_4(1)A}). Using the saddle point method in (\ref{S4_th2_4(1)pnt2}), we obtain (\ref{S2_th2_4}) with $\Lambda(n,t)$ and $\Phi(n,t)=\sqrt{n}\,\Phi_0(s_\ast)$ as in (\ref{S2_th2_4(1)lambda}) and (\ref{S2_th2_4(1)phi}). We note that $a\to (3\sqrt{\rho})^{2/3}$ as $A\to 0$ and that $A$ is an increasing function of $a$, so the above result is valid for $a\geq (3\sqrt{\rho})^{2/3}$.

Alternately, on the scale $n=O(t^{2/3})$ we use the representation (\ref{S2_polla}) with the scaling $z=\sqrt{\rho}\,(1-t^{-1/3}w)$ and $v=t^{-1/3}u$ ($w>0,u>0$). Then we have 
\begin{equation*}
z^{-n}\sim\rho^{-n/2}\,\exp\Big(a\,w\,t^{1/3}+\frac{a\,w^2}{2}\Big),
\end{equation*}
\begin{equation*}
\big[1+\rho-2\sqrt{\rho}\,\cos(v)\big]\,t\sim(1-\sqrt{\rho})^2\,t+\sqrt{\rho}\,u^2\,t^{1/3},
\end{equation*}
\begin{equation*}
\frac{\sin(v)}{1+e^{\pi\cot(v)}}\sim u\,t^{-1/3}\,\exp\Big(-\frac{\pi\,t^{1/3}}{u}\Big),
\end{equation*}
\begin{equation*}
\frac{(1-\sqrt{\rho}\,e^{iv})^{m_0}}{(1-\sqrt{\rho}\,e^{-iv})^{m_0+1}}\sim\frac{1}{1-\sqrt{\rho}}\exp\Big(\frac{\sqrt{\rho}}{1-\sqrt{\rho}}\Big),
\end{equation*}
and
\begin{equation*}
\frac{(\sqrt{\rho}\,e^{-iv}-z)^{m_0}}{(\sqrt{\rho}\,e^{iv}-z)^{m_0+1}}\sim\frac{t^{1/3}}{\sqrt{\rho(w^2+u^2)}}\exp\Big\{\frac{i}{2u}t^{1/3}\log\Big(\frac{w-iu}{w+iu}\Big)+\frac{u^2}{2(w^2+u^2)}\Big\}.
\end{equation*}
It follows that
\begin{eqnarray}\label{S4_th2_4(3)}
p_n(t)&\sim&\frac{2\rho^{-n/2}}{1-\sqrt{\rho}}\exp\Big(\frac{\sqrt{\rho}}{1-\sqrt{\rho}}\Big)\,t^{-2/3}\,e^{-(1-\sqrt{\rho})^2\,t}\frac{1}{2\pi i}\int_{Br_\ast}e^{a\,w\,t^{1/3}+a\,w^2/2}\nonumber\\
&&\times\Big(\int_0^\infty g(u,w)\,e^{-t^{1/3}f(u,w)}\,du\Big)dw,
\end{eqnarray}
where
\begin{equation*}
g(u,w)=\frac{u}{\sqrt{w^2+u^2}}\,\exp\Big\{\frac{u^2}{2(w^2+u^2)}\Big\},
\end{equation*}
\begin{equation*}
f(u,w)=\sqrt{\rho}\,u^2+\frac{\pi}{u}-\frac{i}{2u}\log\Big(\frac{w-iu}{w+iu}\Big),
\end{equation*}
and the contour $Br_\ast$ is a vertical contour in the complex $w$-plane with $\Re(w)$ sufficiently large.
The function $f(u,w)$ has its maximum at $u_\ast=u_\ast(w)$, which satisfies
\begin{equation}\label{S4_th2_4(3)f}
\frac{\partial}{\partial u}f(u,w)=2\sqrt{\rho}\,u-\frac{\pi}{u^2}+\frac{i}{2u^2}\log\Big(\frac{w-iu}{w+iu}\Big)-\frac{w}{u\,(w^2+u^2)}=0.
\end{equation}
Then using the Laplace method in the inner integral of (\ref{S4_th2_4(3)}) implies that
\begin{eqnarray}\label{S4_th2_4(3)2}
p_n(t)&\sim&\frac{2\sqrt{2\pi}}{1-\sqrt{\rho}}\,\rho^{-n/2}\exp\Big(\frac{\sqrt{\rho}}{1-\sqrt{\rho}}\Big)\,t^{-5/6}\,e^{-(1-\sqrt{\rho})^2t}\nonumber\\
&&\times\,\frac{1}{2\pi i}\int_{Br_\ast}\frac{g(u_\ast,w)\,e^{a\,w^2/2}}{\sqrt{f_{uu}(u_\ast,w)}}e^{t^{1/3}[a\,w-f(u_\ast,w)]}dw.
\end{eqnarray}
Let $h(w)=a\,w-f(u_\ast(w),w)$. Then the saddle point equation $h'(w)=0$ along with (\ref{S4_th2_4(3)f}) leads to
\begin{equation}\label{S4_th2_4(3)wu}
w_0^2+u_0^2=\frac{1}{a}.
\end{equation}
Here we denote the solution to $h'(w)=0$ by $w_0=w_0(a)$ and set $u_0=u_\ast(w_0)$. Applying the saddle point method to (\ref{S4_th2_4(3)2}), we find that
\begin{equation}\label{S4_th2_4(3)3}
p_n(t)\sim\rho^{-n/2}\,\frac{2\,t^{-5/6}\,g(u_0,w_0)\,e^{a\,w_0^2/2}}{(1-\sqrt{\rho})\,\sqrt{f_{uu}(u_0,w_0)\,h''(w_0)}}\,e^{-(1-\sqrt{\rho})^2\,t+t^{1/3}\,h(w_0)}.
\end{equation}
If we let $C=1-a\,w_0^2\;(0<C<1)$, then from (\ref{S4_th2_4(3)wu}) it follows that
\begin{equation}\label{S4_th2_4(3)C}
w_0=\sqrt{(1-C)/a}>0\quad \textrm{and} \quad u_0=\sqrt{C/a}>0.
\end{equation}
Since $w_0$ and $u_0$ satisfy (\ref{S4_th2_4(3)f}), using (\ref{S4_th2_4(3)C}) in (\ref{S4_th2_4(3)f}) leads to (\ref{S2_th2_4(3)C}). Hence, after some simplification in (\ref{S4_th2_4(3)3}), we obtain (\ref{S2_th2_4}) with $\Lambda$ and $\Phi=t^{1/3}\,h(w_0)$ as in (\ref{S2_th2_4(3)lambda}) and (\ref{S2_th2_4(3)phi}). We note that $a\to 0$ as $C\to 0$, $a\to (4\sqrt{\rho}/\pi)^{2/3}$ as $C\to 1$ and $C=C(a)$ is an increasing function of $a$. Thus the above result is valid for $0<a\leq (4\sqrt{\rho}/\pi)^{2/3}$.

We now consider the range $(4\sqrt{\rho}/\pi)^{2/3}<a<(3\sqrt{\rho})^{2/3}$. This is difficult to treat using either of the representations in (\ref{S2_th1_phat}) and (\ref{S2_polla}), as the various saddle points become complex. However, we now show that the results in case 4(b) of Theorem 2.2 can be obtained by smoothly continuing the results for case 4(a), or those of case 4(c).

First we consider (\ref{S2_th2_4(1)A}) which we rewrite as 
\begin{eqnarray}\label{S4_4(2)A}
2\sqrt{\rho}\,A\,a^{-3/2}&=&\sqrt{1+A}-\frac{\mathrm{arcsinh}\big(\sqrt{A}\big)}{\sqrt{A}}\\
&=&\frac{2}{3}\,A+O(A^2),\nonumber
\end{eqnarray}
where the right side of (\ref{S4_4(2)A}) is an analytic function of $A$. The curve $a=(3\sqrt{\rho})^{2/3}$ corresponds to $A=0$. Setting $A=-B$ in (\ref{S4_4(2)A}) we obtain (\ref{S2_th2_4(2)B}), which is the analytic continuation of (\ref{S4_4(2)A}) into the range $A<0$. Then (\ref{S2_th2_4(2)lambda}) and (\ref{S2_th2_4(2)phi}) follow by replacing $A$ by $-B$ in (\ref{S2_th2_4(1)lambda}) and (\ref{S2_th2_4(1)phi}). We now show that case 4(b) also follows by the continuation of case 4(c), as $a$ increases past $(4\sqrt{\rho}/\pi)^{2/3}$, which corresponds to $C=1$. The smooth continuation of (\ref{S2_th2_4(3)C}) as $C$ increases past $C=1$ follows by replacing $C$ by $B$, $\arcsin(\sqrt{C})$ by $\pi-\arcsin(\sqrt{B})$ and $\sqrt{1-C}$ by $-\sqrt{1-B}$. Note that viewing $C$ as a function of $(n,t)$, $1-C$ has a double zero along $n\,t^{-2/3}=(4\sqrt{\rho}/\pi)^{2/3}$. These observations show that the three cases in Theorem 2.2 for $n=O(t^{2/3})$ really correspond to a single asymptotic scale. A geometric interpretation of these three cases is given in Section 6.

In the matching region between cases 3 and 4(a) in Theorem 2.2, we let $n/t\to 0$ in (\ref{S2_th2_3}), which yields (\ref{S2_th2_matching34}). On the other hand, letting $a\to\infty$ in (\ref{S2_th2_4(1)A}), we obtain
\begin{equation}\label{S4_matching34}
A=\frac{1}{4\rho}\,a^3-3\log(a)+\big(1+\log(\rho)\big)+O\Big(\frac{1}{a}\Big).
\end{equation}
Using (\ref{S4_matching34}) in (\ref{S2_th2_4(1)lambda}) and (\ref{S2_th2_4(1)phi}) also leads to (\ref{S2_th2_matching34}).

Finally, we consider $n=O(1)$ and $t\to\infty$. We use the representation (\ref{S2_polla}) and scale $v=t^{-1/3}\,V$. Then the inner integral in (\ref{S2_polla}) becomes
\begin{equation*}\label{S4_th2_5}
\frac{2\sqrt{\rho}\,\exp\Big(\frac{\sqrt{\rho}}{1-\sqrt{\rho}}+\frac{\sqrt{\rho}}{\sqrt{\rho}-z}\Big)}{(1-\sqrt{\rho})(\sqrt{\rho}-z)}\,e^{-(1-\sqrt{\rho})^2t}\,t^{-1/3}\int_0^\infty V\,e^{-t^{1/3}(\pi/V+\sqrt{\rho}\,V^2)}\,dV.
\end{equation*}
Using the Laplace method we find that the integrand is maximal at $V=\big(\frac{\pi}{2\sqrt{\rho}}\big)^{1/3}$, and then making the transformation $z\to\sqrt{\rho}\,z$ in the outer integral in (\ref{S2_polla}) lead to (\ref{S2_th2_5}). 

To verify the asymptotic matching between cases 4(c) and 5 in Theorem 2.2, we let $a\to 0$ in (\ref{S2_th2_4(3)C}). It follows that
\begin{equation}\label{S4_matching45_C}
C\sim\Big(\frac{\pi}{2\sqrt{\rho}}\Big)^{2/3}\,a.
\end{equation}
Using (\ref{S4_matching45_C}) in (\ref{S2_th2_4(3)lambda}) and (\ref{S2_th2_4(3)phi}) yields (\ref{S2_th2_matching45}). On the other hand, we can let $n\to\infty$ in (\ref{S2_th2_5}). We scale $z=1-w/\sqrt{n}$ in the contour integral in (\ref{S2_th2_5}) and use the saddle point method. There is a saddle point at $w=1$ and we obtain
\begin{equation}\label{S4_matching45_int}
\frac{1}{2\pi i}\oint_{\mathscr{C^\ast}}\frac{1}{(1-z)\;z^{n+1}}\exp\Big(\frac{1}{1-z}\Big)dz\sim\frac{\sqrt{e}}{2\sqrt{\pi}\,n^{1/4}}\,e^{2\sqrt{n}}.
\end{equation}
This also leads to (\ref{S2_th2_matching45}), which verifies the matching.

\newpage

\section{Asymptotic results for the case $\rho\approx 1$}

Now we consider the case in which the traffic intensity is close to one. Letting $\rho=1-\epsilon$ with $0<\epsilon\ll 1$, we sketch the main points in deriving Theorem 2.3.

First, we consider $n=O(1)$ and $t=O(1)$. Replacing $\rho$ by $1$ in (\ref{S2_polla}) leads to (\ref{S2_th3_1}). On this scale the solution does not simplify much, but there is little probability mass in heavy traffic on the time scale $t=O(1)$.

Next, we consider $n=O(1)$ but very large time scales $t=\sigma/\epsilon^3=O(\epsilon^{-3})$. We use (\ref{S2_polla}) and scale $v=O(\epsilon)$. Then (\ref{S2_th3_2}) is obtained by making the transformation $z\to \sqrt{\rho}\,z$ in the outer integral and using the Laplace method in the inner integral, where the major contribution comes from the point $u(\sigma)$, which satisfies (\ref{S2_th3_2_u}).

To verify the matching between cases 1 and 2 in Theorem 2.3, we let $\sigma\to 0$ in (\ref{S2_th3_2_u}), which yields
\begin{equation}\label{S5_2_u}
u(\sigma)=\Big(\frac{\pi}{4\sigma}\Big)^{1/3}+\frac{2}{3\pi}+O(\sigma^{1/3}).
\end{equation}
Using (\ref{S5_2_u}) in (\ref{S2_th3_2}) leads to (\ref{S2_th3_matching12}). On the other hand, we can let $t\to\infty$ and scale $v=O(t^{-1/3})$ in (\ref{S2_th3_1}). Then using the Laplace method in the inner integral also yields (\ref{S2_th3_matching12}). This implies that there are no other time scales between $t=O(1)$ and $t=O(\epsilon^{-3})$.

For the case $n=\xi/\epsilon=O(\epsilon^{-1})$ and $t=\tau/\epsilon=O(\epsilon^{-1})$, by similar arguments as in the case when $\rho<1$ and $n=O(t)$ in Section 4, we can easily obtain (\ref{S2_th3_3}). We note that if $\xi\to\infty$ and $\tau\to\infty$ but $\xi/\tau=O(1)$, (\ref{S2_th3_3}) is still valid.

For $t=O(\epsilon^{-1})$, we rewrite (\ref{S2_th3_3}) as
\begin{equation}\label{S5_3_rewrite}
p_n(t)=\epsilon\,P_0(\xi,\tau)+\epsilon^2\,P_1(\xi,\tau)+O(\epsilon^3),
\end{equation}
and remove the condition on $n$ by using (\ref{S5_3_rewrite}) in (\ref{S1_p_remove}) with the scaling $n=\xi/\epsilon$. It follows that
\begin{eqnarray*}
p(t)&\sim&\epsilon\int_0^\infty e^{-\xi}\,P_0(\xi,\tau)d\xi+\epsilon^2\int_0^\infty e^{-\xi}\,\Big[P_1(\xi,\tau)-\frac{\xi}{2}\,P_0(\xi,\tau)\Big]d\xi\\
&=& 2\epsilon\,K_0(2\sqrt{\tau})+\frac{\epsilon^2}{3}\Big[(6-\tau)\,K_0(2\sqrt{\tau})-\sqrt{\tau}\,K_1(2\sqrt{\tau})\Big].
\end{eqnarray*}
Here $K_0$ and $K_1$ are the modified Bessel functions. This recovers the result in Morrison \cite{MO} for the time range $t=O(\epsilon^{-1})$, after we take into account that the results in \cite{MO} are for $\Pr\left[\mathbf{V}>t\right]$.

Now we consider larger space-time scales, with $n=\eta/\epsilon^2=O(\epsilon^{-2})$ and $t=\sigma/\epsilon^3=O(\epsilon^{-3})$. We use a similar method as for $\rho<1$ with the scale $n=O(t^{2/3})$. However, the heavy traffic assumption changes some of the saddle point calculations. We note that (\ref{S4_th2_4(1)pnt}) is still valid in the heavy traffic case, but $\theta_c=-\epsilon^2/4+O(\epsilon^3)$ as $\epsilon\to 0$. Then we scale $\theta=\beta\,\epsilon^2\;(-1/4<\beta<0)$ and notice that
\begin{equation*}
z_\pm=1+\frac{1\pm \sqrt{1+4\beta}}{2}\,\epsilon+O(\epsilon^2).
\end{equation*}
Thus, by scaling $z=1+r\,\epsilon$, the inner integral in (\ref{S4_th2_4(1)pnt}), which is $G_n$, is asymptotically given by
\begin{equation}\label{S5_4_Gnint}
G_n\sim \int_{-\infty}^{\frac{1-\sqrt{1+4\beta}}{2}}J_0(r,\beta)\,e^{\psi(r,\beta)/\epsilon}\,dr,
\end{equation}
where
\begin{eqnarray*}
J_0(r,\beta)&=&\frac{2}{\sqrt{(1-2r)^2-(1+4\beta)}}\bigg(\frac{1-2r-\sqrt{1+4\beta}}{1-2r+\sqrt{1+4\beta}}\bigg)^{-\frac{2\beta+1}{2(4\beta+1)^{3/2}}}\\
&&\times\exp\Big\{\frac{r\,(1+3\beta)+\beta\,(1+2\beta)}{(1+4\beta)(r^2-r-\beta)}-\frac{\eta\,r^2}{2}\Big\},
\end{eqnarray*}
and
\begin{equation*}
\psi(r,\beta)=r\,\eta+\frac{1}{\sqrt{1+4\beta}}\log\bigg(\frac{1-2r-\sqrt{1+4\beta}}{1-2r+\sqrt{1+4\beta}}\bigg).
\end{equation*}
The major contribution to the integral in (\ref{S5_4_Gnint}) comes from $r_\ast=r_\ast(\beta)$, which satisfies $\psi_r(r_\ast,\beta)=0$, so that $r_\ast^2-r_\ast-(\beta+1/\eta)=0$. For $\psi$ to be maximal we need 
$$\psi_{rr}(r_\ast,\beta)=\eta^2\,(2r_\ast-1)<0,$$ 
and this implies that
\begin{equation}\label{S5_4(1)rstar}
r_\ast=\frac{1}{2}\,\Big[1-\sqrt{1+4\beta+4/\eta}\Big]<\frac{1}{2}.
\end{equation}
Using the standard Laplace method in (\ref{S5_4_Gnint}), we obtain 
\begin{equation}\label{S5_4_Gn}
G_n\sim\sqrt{2\pi\epsilon}\,J_1(\beta)\,e^{\psi(r_\ast,\beta)/\epsilon},
\end{equation}
where 
\begin{equation*}
J_1(\beta)\equiv\frac{J_0(r_\ast(\beta),\beta)}{\eta\,(1+4\beta+4/\eta)^{1/4}}.
\end{equation*}
Using (\ref{S5_4_Gn}) in (\ref{S4_th2_4(1)pnt}), and the expansion
\begin{equation*}
\frac{z_-}{1-\rho\,z_-}\Big(\frac{z_+-1}{z_--1}\Big)^\alpha\sim \frac{1}{\epsilon\sqrt{|\beta|}}\bigg(\frac{1+\sqrt{1+4\beta}}{1-\sqrt{1+4\beta}}\bigg)^{\frac{1}{{\scriptstyle{\epsilon}}\,\sqrt{1+4\beta}}-\frac{1+2\beta}{2(1+4\beta)^{3/2}}}\exp\Big(\frac{1+2\beta}{1+4\beta}\Big),
\end{equation*}
we have 
\begin{equation}\label{S5_4_pnt}
p_n(t)\sim\frac{\sqrt{2\pi}\,\epsilon^{3/2}}{2\pi i}\int_{Br''}J_1(\beta)\,J_2(\beta)\,e^{\widetilde{\Psi}(\beta)/\epsilon}d\beta,
\end{equation}
where
\begin{equation*}
J_2(\beta)=\frac{1}{\sqrt{|\beta|}}\bigg(\frac{1+\sqrt{1+4\beta}}{1-\sqrt{1+4\beta}}\bigg)^{-\frac{1+2\beta}{2(1+4\beta)^{3/2}}}\exp\Big(\frac{1+2\beta}{1+4\beta}\Big),
\end{equation*}
\begin{equation*}
\widetilde{\Psi}(\beta)=\beta\,\sigma+\frac{1}{\sqrt{1+4\beta}}\,\log\bigg(\frac{1+\sqrt{1+4\beta}}{1-\sqrt{1+4\beta}}\bigg)+\psi(r_\ast(\beta),\beta),
\end{equation*}
and $Br''$ is a vertical contour in the complex $\beta$-plane. Then $\widetilde{\Psi}'(\beta)=0$ implies that there is a saddle point at $\beta_0=\beta_0(\eta,\sigma)$, which satisfies 
\begin{eqnarray}\label{S5_4(1)rstar2}
&&\frac{\sqrt{1+4\beta_0}}{2\beta_0(1-2r_0)}\Big[2r_0+\sigma\beta_0(1+4\beta_0)(1-2r_0)-(1+4\beta_0)(\eta\beta_0+1)\Big]\nonumber\\
&&\quad\quad=\log\bigg(\frac{1-2r_0-\sqrt{1+4\beta_0}}{1-2r_0+\sqrt{1+4\beta_0}}\bigg)+\log\bigg(\frac{1+\sqrt{1+4\beta_0}}{1-\sqrt{1+4\beta_0}}\bigg).
\end{eqnarray}
Here we denote $r_\ast(\beta_0)$ by $r_0$. Thus, from (\ref{S5_4_pnt}), the saddle point method implies that
\begin{equation*}
p_n(t)\sim\epsilon^2\,\frac{J_1(\beta_0)\,J_2(\beta_0)}{\sqrt{\widetilde{\Psi}''(\beta_0)}}\,e^{\widetilde{\Psi}(\beta_0)/\epsilon}.
\end{equation*}
If we let $\beta_0=\widetilde{A}-1/4\;(0<\widetilde{A}<1/4)$, then after some simplification, we have
\begin{equation}\label{S5_4(1)lambda}
\frac{J_1(\beta_0)\,J_2(\beta_0)}{\sqrt{\widetilde{\Psi}''(\beta_0)}}=\widetilde{\Lambda}(\eta,\sigma)
\end{equation}
and
\begin{equation}\label{S5_4(1)phi}
\widetilde{\Psi}(\beta_0)=\widetilde{\Phi}(\eta,\sigma)+\frac{\eta}{2}-\frac{\sigma}{4},
\end{equation} 
where $\widetilde{\Lambda}$ and $\widetilde{\Phi}$ are given by (\ref{S2_th3_4(1)lambda}) and (\ref{S2_th3_4(1)phi}). Equation (\ref{S2_th3_4(1)A}) is derived by using (\ref{S5_4(1)rstar}) and (\ref{S5_4(1)rstar2}). We note from (\ref{S2_th3_4(1)A}) that $\sigma\uparrow\frac{1}{3}\eta^{3/2}-\frac{8}{3}$ as $\widetilde{A}\downarrow 0$. This implies that $\eta>4$, since $\sigma>0$. Hence (\ref{S2_th3_4(1)lambda}) and (\ref{S2_th3_4(1)phi}) are only valid in the range $\sigma<\frac{1}{3}\eta^{3/2}-\frac{8}{3}$ with $\eta>4$. 

Alternately, on the scale $n=O(\epsilon^{-2})$ and $t=O(\epsilon^{-3})$, we use the representation in (\ref{S2_polla}) and scale $v=\epsilon\,\gamma$ and $z=1-\epsilon\,\omega\;(\gamma>0,\omega>0)$. Then we have
\begin{equation*}
\frac{(1-\sqrt{\rho}\,e^{iv})^{m_0}}{(1-\sqrt{\rho}\,e^{-iv})^{m_0+1}}\sim\frac{2}{\epsilon\sqrt{4\gamma^2+1}}\exp\bigg\{\frac{i}{2\epsilon\gamma}\log\Big(\frac{1-2i\gamma}{1+2i\gamma}\Big)-\frac{4\gamma^2+3}{2(4\gamma^2+1)}\bigg\},
\end{equation*}
\begin{eqnarray*}
\frac{(\sqrt{\rho}\,e^{-iv}-z)^{m_0}}{(\sqrt{\rho}\,e^{iv}-z)^{m_0+1}}&\sim&\frac{2}{\sqrt{4\gamma^2+(2\omega-1)^2}}\exp\Big\{\frac{i}{2\epsilon\gamma}\log\Big(\frac{2\omega-1-2i\gamma}{2\omega-1+2i\gamma}\Big)\Big\}\\
&&\times\exp\bigg\{\frac{4\gamma^2-4\omega+3}{2\big[4\gamma^2+(2\omega-1)^2\big]}\bigg\},
\end{eqnarray*}
and (\ref{S2_polla}) is asymptotically given by
\begin{equation}\label{S5_4(3)pnt2int}
p_n(t)\sim\frac{\epsilon}{2\pi i}\int_{\widetilde{Br}}\exp\big(\frac{\eta\omega}{\epsilon}+\frac{\eta\omega^2}{2}\big)\Big(\int_0^\infty l(\gamma,\omega)\,e^{-k(\gamma,\,\omega)/\epsilon}d\gamma\Big)d\omega,
\end{equation} 
where
\begin{eqnarray*}
l(\gamma,\omega)&=&\frac{8\gamma}{\sqrt{(4\gamma^2+1)\big[4\gamma^2+(2\omega-1)^2\big]}}\\
&&\times\exp\bigg\{\frac{4\gamma^2-4\omega+3}{2\big[4\gamma^2+(2\omega-1)^2\big]}-\frac{4\gamma^2+3}{2(4\gamma^2+1)}-\Big(\frac{1}{8}-\frac{\gamma^2}{2}\Big)\sigma\bigg\},
\end{eqnarray*}
\begin{equation*}
k(\gamma,\omega)=\Big(\gamma^2+\frac{1}{4}\Big)\sigma+\frac{\pi}{\gamma}-\frac{i}{2\gamma}\Big[\log\Big(\frac{2\omega-1-2i\gamma}{2\omega-1+2i\gamma}\Big)+\log\Big(\frac{1-2i\gamma}{1+2i\gamma}\Big)\Big],
\end{equation*}
and $\widetilde{Br}$ is a vertical contour in the complex $\omega$-plane with $\Re(\omega)$ sufficiently large.
Thus, the major contribution to the inner integral in (\ref{S5_4(3)pnt2int}) comes from $\gamma_\ast=\gamma_\ast(\omega)$, which satisfies $k_{\gamma}(\gamma_\ast,\omega)=0$, that is 
\begin{eqnarray}\label{S5_4(3)rstar1}
&&2\sigma\gamma_\ast^3-\frac{4\,\omega\,\gamma_\ast\,(4\gamma_\ast^2+2\omega-1)}{(4\gamma_\ast^2+1)\big[4\gamma_\ast^2+(2\omega-1)^2\big]}-\pi\nonumber\\
&&\quad\quad+\frac{i}{2}\Big[\log\Big(\frac{2\omega-1-2i\gamma_\ast}{2\omega-1+2i\gamma_\ast}\Big)+\log\Big(\frac{1-2i\gamma_\ast}{1+2i\gamma_\ast}\Big)\Big]=0.
\end{eqnarray}
Using the Laplace method in the inner integral in (\ref{S5_4(3)pnt2int}) yields
\begin{equation}\label{S5_4(3)pnt}
p_n(t)\sim\frac{\sqrt{2\pi}\,\epsilon^{3/2}}{2\pi i}\int_{\widetilde{Br}}\frac{l(\gamma_\ast,\omega)}{\sqrt{k_{\gamma\gamma}(\gamma_\ast,\omega)}}\,e^{\eta\,\omega^2/2}\,e^{\Phi_1(\omega)/\epsilon}d\omega,
\end{equation}
where $\Phi_1(\omega)=\eta\,\omega-k(\gamma_\ast(\omega),\omega)$. The saddle point equation $\Phi_1'(\omega)=0$ has the solution $\omega_0=\omega_0(\eta,\sigma)$, which satisfies
\begin{equation}\label{S5_4(3)r0}
4\gamma_0^2+(2\omega_0-1)^2=\frac{4}{\eta},\quad \gamma_0=\gamma_\ast(\omega_0).
\end{equation}
In order that $\Phi_1''(\omega_0)>0$, from (\ref{S5_4(3)r0}) it follows that 
\begin{equation}\label{S5_4(3)w0}
\omega_0=\frac{1}{2}+\sqrt{\frac{1}{\eta}-\gamma_0^2}\,.
\end{equation}
Then using the saddle point method in (\ref{S5_4(3)pnt}), we have 
\begin{equation}\label{S5_4(3)pnt2}
p_n(t)\sim\epsilon^2\,\frac{l(\gamma_0,\omega_0)}{\sqrt{k_{\gamma\gamma}(\gamma_0,\omega_0)\,\Phi_1''(\omega_0)}}\,e^{\eta\,\omega_0^2/2}\,e^{\Phi_1(\omega_0)/\epsilon}.
\end{equation}
If let $\gamma_0=\sqrt{\widetilde{C}}\;(0<\widetilde{C}<1/\eta)$, then from (\ref{S5_4(3)w0}) we obtain $\omega_0=1/2+\sqrt{1/\eta-\widetilde{C}}$. Using this in (\ref{S5_4(3)rstar1}) leads to (\ref{S2_th3_4(3)C}). It follows that 
\begin{equation*}
\widetilde{\Lambda}(\eta,\sigma)=\frac{l(\gamma_0,\omega_0)}{\sqrt{k_{\gamma\gamma}(\gamma_0,\omega_0)\,\Phi_1''(\omega_0)}}\,e^{\eta\,\omega_0^2/2}
\end{equation*}
and
\begin{equation*}
\widetilde{\Phi}(\eta,\sigma)=\Phi_1(\omega_0)-\frac{\eta}{2}+\frac{\sigma}{4}
\end{equation*}
in (\ref{S2_th3_4(3)lambda}) and (\ref{S2_th3_4(3)phi}). We note from (\ref{S2_th3_4(3)C}) that at $\widetilde{C}=1/\eta$, we have
\begin{equation}\label{S5_4_curve2}
\sigma=\frac{1}{2}\,\eta^{3/2}\,\Big[\frac{\pi}{2}+\frac{4\sqrt{\eta}}{4+\sqrt{\eta}}-\arcsin\Big(\sqrt{\frac{4}{4+\eta}}\Big)\Big].
\end{equation}  
We also have $\sigma\to \pi/(2\widetilde{C}^{3/2})$ as $\widetilde{C}\to 0$. Hence (\ref{S2_th3_4(3)lambda}) and (\ref{S2_th3_4(3)phi}) are only valid when $\sigma$ exceeds the right side of (\ref{S5_4_curve2}). 

The range of $\sigma$ between
\begin{equation}\label{S5_4_curve1}
\sigma=\frac{1}{3}\,\eta^{3/2}-\frac{8}{3}
\end{equation}
and (\ref{S5_4_curve2}) is difficult to treat using either of the representations in (\ref{S2_th1_phat}) and (\ref{S2_polla}), as the various saddle points become complex. Similarly as in case 4(b) in Theorem 2.2, we now show that the results in case 4(b) of Theorem 2.3 can be obtained by smoothly continuing the results for case 4(a), or those of case 4(c).

First we consider (\ref{S2_th3_4(1)A}) which we rewrite as 
\begin{eqnarray}\label{S5_4(1)A}
2\widetilde{A}\,\sigma&=& -\frac{2}{1-4\widetilde{A}}-\frac{1}{{\sqrt{\widetilde{A}}}}\,\mathrm{arcsinh}\Big(\sqrt{\widetilde{A}\,\eta}\Big)+\frac{1}{{\sqrt{\widetilde{A}}}}\,\mathrm{arcsinh}\bigg(\sqrt{\frac{4\widetilde{A}}{1-4\widetilde{A}}\bigg)}\nonumber\\
&&+\;\sqrt{\eta\,(1+\widetilde{A}\,\eta)}\\
&=& \Big(\frac{2}{3}\,\eta^{3/2}-\frac{16}{3}\Big)\,\widetilde{A}+O(\widetilde{A}^2),\nonumber
\end{eqnarray}
where the right side of (\ref{S5_4(1)A}) is an analytic function of $\widetilde{A}$. The curve (\ref{S5_4_curve1}) corresponds to $\widetilde{A}=0$. Setting $\widetilde{A}=-\widetilde{B}$ in (\ref{S5_4(1)A}) we obtain (\ref{S2_th3_4(2)B}), which is the analytic continuation of (\ref{S5_4(1)A}) into the range $\widetilde{A}<0$. Then (\ref{S2_th2_4(3)lambda}) and (\ref{S2_th2_4(3)phi}) follow by replacing $\widetilde{A}$ by $-\widetilde{B}$ in (\ref{S2_th3_4(1)lambda}) and (\ref{S2_th3_4(1)phi}). We now show that case 4(b) also follows by the continuation of case 4(c), as $\sigma$ decreases past the curve (\ref{S5_4_curve2}), which corresponds to $\widetilde{C}= 1/\eta$. The smooth continuation of (\ref{S2_th3_4(3)C}) as $\widetilde{C}$ increases past $\widetilde{C}= 1/\eta$ follows by replacing $\widetilde{C}$ by $\widetilde{B}$, $\arcsin(\sqrt{\widetilde{C}\,\eta})$ by $\pi-\arcsin(\sqrt{\widetilde{B}\,\eta})$ and $\sqrt{1-\widetilde{C}\,\eta}$ by $-\sqrt{1-\widetilde{B}\,\eta}$. Note that viewing $\widetilde{C}$ as a function of $(\eta,\sigma)$, $1-\widetilde{C}\,\eta$ has a double zero along the curve (\ref{S5_4_curve2}). These observations show that the three cases in item 4 of Theorem 2.3 really correspond to a single asymptotic scale. A geometric interpretation of these three cases is also given in Section 6.

Now we consider the matching between cases 3 and 4(a) in Theorem 2.3. If we fix $\sigma$ but let $\eta=\zeta/\epsilon\to\infty$ in (\ref{S2_th3_4(1)A}), it follows that
\begin{equation}\label{S5_4(1)match}
\widetilde{A}=\frac{1}{4}-\frac{1}{\zeta}\,\epsilon-\frac{\sigma}{\zeta^2}\,\epsilon^2+O(\epsilon^3).
\end{equation}
Using (\ref{S5_4(1)match}) in (\ref{S2_th3_4(1)lambda}) and (\ref{S2_th3_4(1)phi}), we have 
\begin{equation*}
p_n(t)\sim\frac{\epsilon^3}{\zeta}e^{-\sigma/\zeta}=\frac{\epsilon}{\xi}e^{-\tau/\xi}.
\end{equation*}
Alternately, if we fix $\eta$ but let $\sigma=\tau\epsilon^2\to 0$ in (\ref{S2_th3_4(1)A}), we have
\begin{equation}\label{S5_4(1)match2}
\widetilde{A}=\Big(\frac{1}{4}-\frac{1}{\eta}\Big)-\frac{\tau}{\eta^2}\,\epsilon^2+O(\epsilon^4).
\end{equation}
Then $p_n(t)$ becomes  
\begin{equation*}
p_n(t)\sim\frac{\epsilon^2}{\eta}e^{-\tau\epsilon/\eta}=\frac{\epsilon}{\xi}e^{-\tau/\xi}.
\end{equation*}
These calculations verify the matching. 

Next, we consider the matching between cases 2 and 4(c) in Theorem 2.3. Since $\gamma_0$ and $\omega_0$ satisfy (\ref{S5_4(3)rstar1}), we use (\ref{S5_4(3)w0}) in (\ref{S5_4(3)rstar1}) and let $\eta\to 0$. Then the leading term in the asymptotic expansion in (\ref{S5_4(3)rstar1}) leads to (\ref{S2_th3_2_u}) with $u=\gamma_0$. Then letting $\eta\to 0$ in (\ref{S2_th3_4(3)lambda}) and (\ref{S2_th3_4(3)phi}) leads to (\ref{S2_th3_matching24}). On the other hand, if we scale $z=1-\epsilon\,w$ in (\ref{S2_th3_2}) and use the saddle point method, we obtain
 \begin{equation}\label{S5_match24}
\frac{1}{2\pi i}\oint_{\mathscr{C^\ast}}\frac{1}{(1-z)\;z^{n+1}}\exp\Big(\frac{1}{1-z}\Big)dz\sim\frac{\sqrt{\epsilon}}{2\sqrt{\pi}}\,\eta^{-1/4}\exp\Big(\frac{2\sqrt{\eta}}{\epsilon}+\frac{1}{2}\Big).
\end{equation}
Using (\ref{S5_match24}) in (\ref{S2_th3_2}) also leads to (\ref{S2_th3_matching24}). This verifies the matching.

To remove the condition on $n$, we use (\ref{S5_4(3)pnt2}) in (\ref{S1_p_remove}). Since $\rho^n\sim\exp(\eta/\epsilon-\eta/2)$, it follows that
\begin{equation}\label{S5_pt}
p(t)\sim\epsilon\int_0^\infty\frac{l(\gamma_0,\omega_0)}{\sqrt{k_{\gamma\gamma}(\gamma_0,\omega_0)\,\Phi_1''(\omega_0)}}\exp\Big\{-\frac{1}{\epsilon}\Omega(\eta,\sigma)+\Big(\frac{\eta\,\omega_0^2}{2}-\frac{\eta}{2}\Big)\Big\}d\eta.
\end{equation}
Here $\Omega(\eta,\sigma)=\eta-\Phi_1(\omega_0(\eta,\sigma))$. Then $\Omega_\eta(\eta,\sigma)=0$ implies that the major contribution comes from $\eta=\eta_0(\sigma)$ which satisfies $\omega_0(\eta_0,\sigma)=1$. Then using the Laplace method in (\ref{S5_pt}), we have
\begin{equation}\label{S5_ptasym}
p(t)\sim\epsilon^{3/2}\,\frac{\sqrt{2\pi}\;l(\gamma_0,\omega_0)}{\sqrt{k_{\gamma\gamma}(\gamma_0,\omega_0)\,\Phi_1''(\omega_0)\,\Omega_{\eta\eta}(\eta_0,\sigma)}}\,e^{-\Omega(\eta_0,\,\sigma)/\epsilon}.
\end{equation}
Here we set $\omega_0\equiv\omega_0(\eta_0,\sigma)$ and $\gamma_0\equiv\gamma_\ast(\omega_0(\eta_0,\sigma))$. From (\ref{S5_4(3)w0}), if we let $\gamma_0=\frac{1}{2}\cot\big(\frac{\psi}{2}\big)\;(0<\psi<\pi)$, then $\eta_0=4\sin^2\big(\frac{\psi}{2}\big)$. Thus, (\ref{S5_4(3)rstar1}) becomes 
\begin{equation*}
\sigma=\sigma(\psi)=4\big(\sin(\psi)+\psi\big)\,\tan^3\Big(\frac{\psi}{2}\Big)
\end{equation*}
and we also have 
\begin{equation*}
\Omega(\eta_0,\sigma)=-2\psi\,\tan\Big(\frac{\psi}{2}\Big)+\frac{\sigma}{4}\,\csc^2\Big(\frac{\psi}{2}\Big)\equiv F_0(\psi).
\end{equation*}
Letting 
\begin{equation*}
F_1(\psi)\equiv\frac{\sigma}{8}\csc^2\Big(\frac{\psi}{2}\Big)-\psi\,\tan\Big(\frac{\psi}{2}\Big),
\end{equation*}
(\ref{S5_ptasym}) becomes 
\begin{equation*}
p(t)\sim\epsilon^{3/2}\,\cot\Big(\frac{\psi}{2}\Big)\,\sqrt{\frac{2\pi}{F_0''(\psi)}}\,e^{-F_1(\psi)}\,e^{-F_0(\psi)/\epsilon},
\end{equation*}
which recovers the result in Morrison \cite{MO} on the time scale $t=O(\epsilon^{-3})$.

\newpage

\section{Singular perturbation method}

Now we discuss an alternate, singular perturbation approach for deriving the asymptotic approximations.

We first assume that the traffic intensity $\rho$ is fixed and less than one. We introduce a small parameter $\delta$ ($0<\delta\ll 1$), let $n=N/\delta$, $t=T/\delta$, and expand $p_n(t)$ as follows:
\begin{equation}\label{S6_1_pexpansion}
p_n(t)=\delta q_0(N,T)+\delta^2 q_1(N,T)+O(\delta^3).
\end{equation}
Using the recurrence equation (\ref{S2_recu1}), the leading term in (\ref{S6_1_pexpansion}) satisfies
\begin{equation*}\label{S6_1_pdeleading}
\frac{\partial q_0}{\partial T}=-(1-\rho)\frac{\partial q_0}{\partial N}-\frac{1}{N}q_0
\end{equation*}
with the initial condition $q_0(N,0)=1/N$. Solving this PDE by the method of characteristics, we obtain
\begin{equation}\label{S6_1_q0}
q_0(N,T)=\frac{1}{N}\Big[1-(1-\rho)\frac{T}{N}\Big]^{\frac{\rho}{1-\rho}}=\frac{1}{N}\Delta_1^{\frac{\rho}{1-\rho}},\quad \Delta_1>0.
\end{equation}
Using (\ref{S6_1_q0}), we can solve for $q_1=q_1(N,T)$, which satisfies the PDE
\begin{equation*}
\frac{\partial q_1}{\partial T}=\frac{1+\rho}{2}\frac{\partial^2 q_0}{\partial^2 N}+\frac{1}{N}\frac{\partial q_0}{\partial N}-(1-\rho)\frac{\partial q_1}{\partial N}+\frac{1}{N^2}q_0-\frac{1}{N}q_1,
\end{equation*}
with the initial conditon $q_1(N,0)=-1/N^2$.
Hence, we have
\begin{eqnarray}\label{S6_1_q1}
q_1(N,T)&=&\frac{1}{2(1-\rho)^3N^2}\,\Delta_1^{\frac{3\rho-2}{1-\rho}}\,\Big[\rho(2\rho^2+\rho-1)+4\rho^2\Delta_1\log(\Delta_1)\nonumber\\
&&+\,6\rho(1-\rho)\Delta_1-(\rho^2-\rho+2)\Delta_1^2\Big].
\end{eqnarray}
Using (\ref{S6_1_q0}) and (\ref{S6_1_q1}) in (\ref{S6_1_pexpansion}), we get (\ref{S2_th2_1}) upon setting $\delta=1$, so that $(N,T)=(n,t)$. From (\ref{S6_1_q0}) we note that this result is valid for $n,\,t\to\infty$ with $n/t>1-\rho$.

Next, we consider $n,\,t\to\infty$ with $0<n/t<1-\rho$. We assume that $p_n(t)$ has an expansion in the following form
\begin{equation*}
p_n(t)=\delta^{\nu}\,e^{\Theta(N,T)/\delta}\,\big[K^{(1)}(N,T)+\delta\,K^{(2)}(N,T)+O(\delta^2)\big].
\end{equation*}
Using this in (\ref{S2_recu1}) yields the PDEs
\begin{equation}\label{S6_3_Thetapde}
\Theta_T=\rho\,e^{\Theta_N}+e^{-\Theta_N}-1-\rho,
\end{equation}
and
\begin{eqnarray}\label{S6_3_Kpde}
K^{(1)}_T&=&(\rho\,e^{\Theta_N}-e^{-\Theta_N})\,K^{(1)}_N\nonumber\\
&&+\Big(\frac{\rho}{2}\,\Theta_{NN}\,e^{\Theta_N}+\frac{1}{2}\,\Theta_{NN}\,e^{-\Theta_N}-\frac{1}{N}\,e^{-\Theta_N}\Big)\,K^{(1)}.
\end{eqnarray}
The PDE (\ref{S6_3_Thetapde}) can be solved by the method of characteristics with all of the rays starting from the origin $(N,T)=(0,0)$. This leads to   
\begin{equation}\label{S6_3_Theta}
\Theta(N,T)=T\Big(-1-\rho+\sqrt{\frac{N^2}{T^2}+4\rho}\Big)+N\log\Big[\frac{1}{2\rho}\Big(-\frac{N}{T}+\sqrt{\frac{N^2}{T^2}+4\rho}\Big)\Big]. 
\end{equation}
The function $K^{(1)}(N,T)$ cannot be determined completely, but from (\ref{S6_3_Kpde}) we find that it has the form
\begin{equation}\label{S6_3_K}
K^{(1)}(N,T)=N^{-1-\frac{1}{2}\sqrt{1+4\rho\,T^2/N^2}}\,K\Big(\frac{N}{T}\Big).
\end{equation}
Since $p_n(t)$ must ultimately be independent of $\delta$, we can set $\delta=1$. Thus, we have
\begin{equation}\label{S6_3_pnt}
p_n(t)\sim n^{-1-\frac{1}{2}\sqrt{1+4\rho\,t^2/n^2}}\,K\Big(\frac{n}{t}\Big)\,e^{\Theta(n,t)},\quad 0<\frac{n}{t}<1-\rho.
\end{equation}

On the scale $n,t\to\infty$ with $n/t=1-\rho+O(t^{-1/2})$, we let $N=(1-\rho)\,T+\sqrt{\delta}\,S$. From (\ref{S6_1_q0}), it follows that
\begin{equation*}
\frac{1}{N}\Big[1 - (1 - \rho )\frac{T}{N}\Big]^{\frac{\rho }{{1 - \rho }}}  = \delta ^{\frac{{\rho }}{{2(1 - \rho )}}} \,S^{\frac{\rho }{{1 - \rho }}} \,N^{ - \frac{1}{{1 - \rho }}}. 
\end{equation*}
Then we expand $p_n(t)$ in the form
\begin{equation}\label{S6_2_pnt}
p_n(t)\sim \delta^{\nu_0}\,N^{-\frac{1}{1-\rho}}\,\mathscr{F}(N,S).
\end{equation}
Using (\ref{S6_2_pnt}) in (\ref{S2_recu1}), we obtain for $\mathscr{F}$ the heat equation
\begin{equation}\label{S6_2_Fpde}
\frac{1+\rho}{2}\,\mathscr{F}_{SS}=(1-\rho)\,\mathscr{F}_N.
\end{equation}
We next obtain two matching conditions between the scale $S=O(1)$ and the ranges $n/t>1-\rho$ and $0<n/t<1-\rho$. We first let $N=(1-\rho)T+\sqrt{\delta}\,S$ in (\ref{S6_1_pexpansion}) and (\ref{S6_1_q0}), which yields
\begin{equation*}
\delta\,q_0(N,T)=\delta^{\frac{2-\rho}{2(1-\rho)}}\,N^{-\frac{1}{1-\rho}}\,S^{\frac{\rho}{1-\rho}}.
\end{equation*}
If (\ref{S6_1_pexpansion}) and (\ref{S6_2_pnt}) were to match, the above should agree with the behavior of (\ref{S6_2_pnt}) as $S\to\infty$, which implies that $\nu_0=\frac{2-\rho}{2(1-\rho)}$ and
\begin{equation}\label{S6_2_match_sinfty}
\mathscr{F}(N,S)\,\sim\,S^{\frac{\rho}{1-\rho}}\quad \textrm{as}\quad S\to \infty.
\end{equation}

We next consider the matching between the scales $n/t\approx 1-\rho$ and $0<n/t<1-\rho$. We let $n=N/\delta$ and $t=T/\delta$ in (\ref{S6_3_pnt}) and let $N/T\to1-\rho$. In this limit we have
\begin{equation*}
\Theta(n,t)=\frac{1}{\delta}\,\Theta(N,T)\sim-\frac{1}{\delta}\,\frac{(1-\rho)}{2(1+\rho)}\Big[N-(1-\rho)\,T\Big]^2,
\end{equation*}
so that
\begin{equation*}
p_n(t)\sim\Big(\frac{N}{\delta}\Big)^{-\frac{3-\rho}{2(1-\rho)}}\,\exp\Big\{-\frac{(1-\rho)\,S^2}{2(1+\rho)\,N}\Big\}\,K\Big(\frac{N}{T}\Big)\Big|_{N/T\to 1-\rho}.
\end{equation*}
We furthermore assume that $K(\cdot)$ has some algebraic behavior as $N/T\to 1-\rho$, in the form 
\begin{eqnarray*}
K\Big(\frac{N}{T}\Big)&\sim& c_0\,\Big[(1-\rho)-\frac{N}{T}\Big]^{\nu_1}\\
&=& c_0\,\Big[-(1-\rho)\sqrt{\delta}\frac{S}{N}\Big]^{\nu_1}\quad \textrm{as}\quad \frac{N}{T}\to 1-\rho,
\end{eqnarray*}
where $c_0$ and $\nu_1$ are constants that will be determined later.
We thus obtain the second matching condition
\begin{eqnarray}\label{S6_2_matching}
&&\delta^{\nu_0}\,N^{-\frac{1}{1-\rho}}\,\mathscr{F}(N,S)\nonumber\\
&&\sim\,c_0\,\Big(\frac{N}{\delta}\Big)^{-\frac{3-\rho}{2(1-\rho)}}\,\Big[-(1-\rho)\sqrt{\delta}\frac{S}{N}\Big]^{\nu_1}\,\exp\Big\{-\frac{(1-\rho)\,S^2}{2(1+\rho)\,N}\Big\},
\end{eqnarray}
as $S\to-\infty$. By comparing powers of $\delta$ in (\ref{S6_2_matching}), it follows that 
$$ \nu_0=\frac{3-\rho}{2(1-\rho)}+\frac{\nu_1}{2}$$
and thus $\nu_1=-\frac{1}{1-\rho}$.
The matching conditions also suggest that we seek a solution of (\ref{S6_2_Fpde}) in terms of the similarity variable $S/\sqrt{N}$, with
\begin{equation}\label{S6_2_U}
\mathscr{F}(N,S)=N^{\nu_2}\,\mathscr{H}\Big(\frac{S}{\sqrt{N}}\Big)=N^{\nu_2}\,\mathscr{H}(\Delta_2),
\end{equation}
where 
$$\nu_2=\frac{1}{1-\rho}-\frac{3-\rho}{2(1-\rho)}-\frac{\nu_1}{2}=\frac{\rho}{2(1-\rho)}.$$
Then (\ref{S6_2_match_sinfty}) and (\ref{S6_2_matching}) give the behavior of $\mathscr{H}(\Delta_2)$ as $\Delta_2\to\pm\infty$, as 
\begin{equation*}
\mathscr{H}(\Delta_2)\sim \Delta_2^{\frac{\rho}{1-\rho}}\quad\textrm{as}\quad \Delta_2\to+\infty,
\end{equation*}
\begin{eqnarray}\label{S6_2_matchcon1}
\mathscr{H}(\Delta_2)&\sim& c_0\, (1-\rho)^{-\frac{1}{1-\rho}}\,(-\Delta_2)^{-\frac{1}{1-\rho}}\nonumber\\
&&\times\,\exp\Big\{-\frac{1-\rho}{2(1+\rho)}\Delta_2^2\Big\}\quad\textrm{as}\quad \Delta_2\to-\infty.
\end{eqnarray}
Using (\ref{S6_2_U}) in (\ref{S6_2_Fpde}), we find that $\mathscr{H}$ satisfies the parabolic cylinder equation
\begin{equation*}
(1+\rho)\,\mathscr{H}''(\Delta_2)+(1-\rho)\,\Delta_2\,\mathscr{H}'(\Delta_2)-\rho\,\mathscr{H}(\Delta_2)=0.
\end{equation*}
The solution that satisfies both matching conditions is given by
\begin{equation}\label{S6_2_V}
\mathscr{H}(\Delta_2)=\sqrt{\frac{1-\rho}{1+\rho}}\,\frac{1}{\sqrt{2\pi}}\int_0^{\infty}y^{\frac{\rho}{1-\rho}}\exp\Big\{-\frac{1-\rho}{2(1+\rho)}(y-\Delta_2)^2\Big\}dy.
\end{equation}
Combining (\ref{S6_2_pnt}), (\ref{S6_2_U}) and (\ref{S6_2_V}) and setting $\delta=1$ we regain (\ref{S2_th2_2}). We also note that by letting $\Delta_2\to-\infty$ in (\ref{S6_2_V}) and using (\ref{S6_2_matchcon1}), we determine $c_0$ as
\begin{equation*}
c_0=\frac{1}{\sqrt{2\pi}}\,\sqrt{\frac{1-\rho}{1+\rho}}\,(1+\rho)^{\frac{1}{1-\rho}}\,\Gamma\Big(\frac{1}{1-\rho}\Big).
\end{equation*}

Now we consider $n,t\to\infty$ and $n\,t^{-2/3}\equiv a=O(1)$. From the case $n,t\to\infty$ with $0<n/t<1-\rho$, if we let $n/t\to 0$ in (\ref{S6_3_Theta}), it follows that
\begin{equation*}
e^{\Theta(n,t)}= \rho^{-n/2}\,\exp\Big\{-(1-\sqrt{\rho})^2\,t+O(n^2/t)\Big\}.
\end{equation*}
Then we set
\begin{equation}\label{S6_4_pnt_form}
p_n(t)=\rho^{-n/2}\,e^{-(1-\sqrt{\rho})^2\,t}\,R_n(t).
\end{equation}
Using (\ref{S6_4_pnt_form}) in (\ref{S2_recu1}), we have the following recurrence equation
\begin{equation}\label{S6_4_Rn_recu}
\frac{1}{\sqrt{\rho}}R_n'(t)=R_{n+1}(t)-2\,R_n(t)+\frac{n}{n+1}\,R_{n-1}(t).
\end{equation}
Letting $n=Y\,\delta^{-2/3}$ and $t=T/\delta$ we assume that $R_n(t)$ has the following asymptotic expansion:
\begin{equation}\label{S6_4_Rn_asym}
R_n(t)\sim \delta^{\nu_3}\,\Delta(Y,T)\,\exp\big\{\delta^{-1/3}\,\Psi(Y,T)\big\}.
\end{equation}
Using (\ref{S6_4_Rn_asym}) in (\ref{S6_4_Rn_recu}), the perturbation method yields the PDEs
\begin{equation}\label{S6_4_pde1}
\frac{1}{\sqrt{\rho}}\Psi_T=(\Psi_Y)^2-\frac{1}{Y}
\end{equation}
and
\begin{equation}\label{S6_4_pdedelta}
\Delta_{T}-2\Psi_Y\,\Delta_Y=\Big(\Psi_{YY}+\frac{1}{Y}\,\Psi_Y\Big)\,\Delta.
\end{equation}

We use the method of characteristics to solve (\ref{S6_4_pde1}), with all of the rays coming from the origin $(Y,T)=(0,0)$. We find that the geometry of the rays naturally defines three regions in the $(Y,T)$ plane (as shown in Figure 1). The first region corresponds to $\Psi_T\geq 0$ and $\Psi_Y<0$, where the rays $A\geq 0$ satisfy (\ref{S2_th2_4(1)A}). Here $\Psi_T=\sqrt{\rho}A/Y$ is constant along a ray. When $A=0$, $a=n\,t^{-2/3}=(3\sqrt{\rho})^{2/3}$, which is the ray denoted by the dashed curve in Figure 1. The second region corresponds to $\Psi_T<0$ and $\Psi_Y<0$, where the rays $0<B<1$ satisfy (\ref{S2_th2_4(2)B}), with now $\Psi_T=-\sqrt{\rho}B/Y$ being constant along a ray. If let $B\to 1$ in (\ref{S2_th2_4(2)B}), then $a\to (4\sqrt{\rho}/\pi)^{2/3}$. $a=(4\sqrt{\rho}/\pi)^{2/3}$ is not a ray, which we denote as a dotted curve in Figure 1. This curve corresponds to the locus of the maximum values of $Y$ achieved along the rays that start from $(0,0)$ and return to $Y=0$ at some later $T>0$. The third region corresponds to $\Psi_T<0$ and $\Psi_Y>0$, where the rays $0<C<1$ satisfy (\ref{S2_th2_4(3)C}). Thus, (\ref{S2_th2_4(1)phi}), (\ref{S2_th2_4(2)phi}) and (\ref{S2_th2_4(3)phi}) are obtained from the corresponding region with $\delta^{-1/3}\Psi(Y,T)=\Phi(n,t)$.

The function $\Delta(Y,T)$ cannot be determined completely from (\ref{S6_4_pdedelta}), but we find that it must have the form
\begin{equation*}
\Delta(Y,T)=\frac{1}{T}\,\Delta_0\Big(\frac{Y}{T^{2/3}}\Big)=\frac{1}{T}\,\Delta_0(a),
\end{equation*}
We shall obtain the behaviors of $\Delta_0(a)$ as $a\to 0$ and $a\to \infty$ below. Thus, by setting $\delta=1$, $p_n(t)$ has the following asymptotic approximation:
\begin{equation}\label{S6_4_pnt}
p_n(t)\sim\frac{\rho^{-n/2}}{t}\Delta_0(a)\;\exp\Big\{-(1-\sqrt{\rho})^2t+\Phi(n,t)\Big\},
\end{equation}
where $\Phi(n,t)$ has three different forms, but passes smoothly through the two transition curves in Figure 1.

For the scale $n=O(1)$ and $t\to\infty$, we assume that $p_n(t)$ has the form
\begin{equation*}
p_n(t)\sim\delta^{\nu_4}\,\rho^{-n/2}\,\mathscr{P}_n(T)\,e^{-(1-\sqrt{\rho})^2\,t}.
\end{equation*}
Using the above expansion in (\ref{S2_recu1}), we obtain the difference equation 
\begin{equation}\label{S6_5_diff}
\mathscr{P}_{n+1}-2\mathscr{P}_n+\frac{n}{n+1}\mathscr{P}_{n-1}=0.
\end{equation}
Using generating functions we can express $\mathscr{P}_n(T)$ in terms of $\mathscr{P}_0(T)$, as
\begin{equation}\label{S6_5_Pn}
\mathscr{P}_n(T)=\mathscr{P}_0(T)\frac{1}{2\pi i}\oint_{\mathscr{C}^\ast}\frac{e^{-1}}{z^{n+1}\,(1-z)}e^{\frac{1}{1-z}}dz.
\end{equation}
Then, by setting $\delta=1$, we have for $n=O(1)$
\begin{equation}\label{S6_5_pnt}
p_n(t)\sim\rho^{-n/2}\,e^{-(1-\sqrt{\rho})^2\,t}\mathscr{P}_0(t)\frac{1}{2\pi i}\oint_{\mathscr{C}^\ast}\frac{e^{-1}}{z^{n+1}\,(1-z)}e^{\frac{1}{1-z}}dz.
\end{equation}

In (\ref{S6_3_pnt}), (\ref{S6_4_pnt}) and (\ref{S6_5_pnt}), the functions $K$, $\Delta_0$ and $\mathscr{P}_0$ cannot be completely determined by the perturbation method. But by examing the matching between the scales, we can find their structure in the matching regions. 

We first consider the matching region between the scales $0<n/t<1-\rho$ and $n=O(t^{2/3})$. If let $n/t\to 0$ in (\ref{S6_3_pnt}), we have 
\begin{equation*}
n^{-1-\frac{1}{2}\sqrt{1+4\rho\,t^2/n^2}}\,e^{\Theta(n,t)}\sim\frac{\rho^{-n/2}}{n}\,e^{-(1-\sqrt{\rho})^2t}\,\exp\Big\{-\frac{n^2}{4\sqrt{\rho}\,t}-\sqrt{\rho}\frac{t}{n}\log(n)\Big\}
\end{equation*}
and we assume that $K(n/t)$ has the form
\begin{equation*}
K\Big(\frac{n}{t}\Big)\sim c_1\Big(\frac{n}{t}\Big)^{\beta_1}\,e^{\gamma_1(n,t)\,t/n}\quad\textrm{as}\quad \frac{n}{t}\to 0.
\end{equation*} 
If let $a=n\,t^{-2/3}\to\infty$ in (\ref{S6_4_pnt}) in the first region where $a>(3\sqrt{\rho})^{2/3}$, we can use (\ref{S4_matching34}) and obtain
\begin{equation*}
\Phi(n,t)\sim-\frac{n^2}{4\sqrt{\rho}\,t}+\sqrt{\rho}\frac{t}{n}\big[\log(\rho)-1+2\log(t)-3\log(n)\big].
\end{equation*}
Then we assume that $\Delta_0(a)$ has an algebraic behavior as $a\to \infty$, in the form
\begin{equation*}
\Delta_0(a)\sim c_2\Big(\frac{n}{t^{2/3}}\Big)^{\beta_2}\quad\textrm{as}\quad a=n\,t^{-2/3}\to \infty.
\end{equation*}
If these two scales are to match in an intermediate limit where $n/t\to 0$ and $a\to\infty$, it follows that
\begin{equation*}
c_1\,\frac{1}{n}\,\Big(\frac{n}{t}\Big)^{\beta_1}=c_2\,\frac{1}{t}\,\Big(\frac{n}{t^{2/3}}\Big)^{\beta_2}.
\end{equation*}
Thus, we conclude that $c_1=c_2$, $\beta_1=1$ and $\beta_2=0$. We also conclude that 
\begin{equation*}
\gamma_1(n,t)=\sqrt{\rho}\big[\log(\rho)-1+2\log(t)-2\log(n)\big].
\end{equation*}

We now consider the  matching region between the scales $n=O(t^{2/3})$ and $n=O(1),\;t\to\infty$. If we let $a\to 0$ in (\ref{S6_4_pnt}) in the third region where $0<a<(4\sqrt{\rho}/\pi)^{2/3}$, we have (\ref{S4_matching45_C}), which implies that
\begin{equation*}
\Phi(n,t)\sim -3\Big(\frac{\pi}{2}\Big)^{2/3}\rho^{1/6}\,t^{1/3}+2\sqrt{n}.
\end{equation*}
We then assume that $\Delta_0(a)$ has the form
\begin{equation*}
\Delta_0(a)\sim c_3\Big(\frac{n}{t^{2/3}}\Big)^{\beta_3}\quad\textrm{as}\quad a=n\,t^{-2/3}\to 0.
\end{equation*}
We let $n\to \infty$ in (\ref{S6_5_pnt}) and use (\ref{S4_matching45_int}), and also assume the following form for $\mathscr{P}_0(t)$ 
\begin{equation*}
\mathscr{P}_0(t)\sim c_4\,t^{\beta_4}\exp\big(\gamma_4\,t^{1/3}\big)\quad\textrm{as}\quad t\to \infty.
\end{equation*}
Note that the exponential factor is indicated by the behavior of $\Phi$ as $a\to 0$.
Then the matching holds provided that 
\begin{equation*}
c_3\,\frac{1}{t}\,\Big(\frac{n}{t^{2/3}}\Big)^{\beta_3}=c_4\,\frac{e^{-1/2}}{2\sqrt{\pi}\,n^{1/4}}\,t^{\beta_4}.
\end{equation*}
Thus, we conclude that 
\begin{equation*}
c_3=\frac{e^{-1/2}}{2\sqrt{\pi}}\,c_4,
\end{equation*}
$\beta_3=-1/4$, $\beta_4=-5/6$ and
\begin{equation*}
\gamma_4=-3\Big(\frac{\pi}{2}\Big)^{2/3}\rho^{1/6}.
\end{equation*}

We have thus shown that much, but certainly not all, of the asymptotic structure of $p_n(t)$ in Theorem 2.2 can be obtained by perturbation methods, which make no recourse to the exact solutions in Theorem 2.1 and (\ref{S2_polla}). We can use this method to obtain the full asymptotic series for the scales for cases 1 and 2 in Theorem 2.2. For the scales in cases 3, 4 and 5 we can obtain partial information only. Specifically, we can get (\ref{S2_th2_3}) only up to the unknown function $K(n/t)$, for which we can infer the behavior as $n/t\to 1-\rho$ and as $n/t\to 0$ (up to the constant $c_1$). Of course this function was fully determined by the saddle point method, as given in (\ref{S2_th2_3_K}). From (\ref{S2_th2_3_K}) we can easily show that the behaviors as $n/t\to 1-\rho$ and as $n/t\to 0$ were correctly predicted by the matching arguments.

Similarly, the perturbation method yielded the function $\Phi(n,t)$ in (\ref{S2_th2_4}) completely, for all 3 ranges of $a$. But, only partial information could be obtained about $\Lambda(n,t)=t^{-1}\Delta_0(a)$. Specifically, we obtained the behavior of $\Delta_0$ as $a\to\infty$ (up to the constant $c_2$) and as $a\to 0$ (up to the constant $c_3$). For $t\to\infty$ with $n=O(1)$ we could determine the expansion of $p_n(t)$ up to the constant $c_4$, which is expressible in terms of $c_3$. For the scale $n=O(t^{2/3})$ the perturbation method led to a nice geometric interpretation of the 3 sub-cases in item 4 of the Theorem 2.2.

Next we very briefly discuss the heavy traffic case via perturbation expansions. We set $\rho=1-\epsilon$ with $0<\epsilon\ll 1$, and first consider $n=\xi/\epsilon$ and $t=\tau/\epsilon$,  expanding $p_n(t)$ as follows: 
\begin{equation}\label{S6_ht_3_P}
p_n(t)=\epsilon\,P_0(\xi,\tau)+\epsilon^2\,P_1(\xi,\tau)+O(\epsilon^3).
\end{equation}
Using (\ref{S6_ht_3_P}) in (\ref{S2_recu1}) we obtain to leading order
\begin{equation}\label{S6_ht_3_P0pde}
\frac{\partial P_0}{\partial\tau}=-\frac{1}{\xi}\,P_0\quad\textrm{with}\quad P_0(\xi,0)=\frac{1}{\xi}
\end{equation}
so that
\begin{equation}\label{S6_ht_3_P0}
P_0(\xi,\tau)=\frac{1}{\xi}\,e^{-\tau/\xi}.
\end{equation}
To obtain $P_1(\xi,\tau)$, we need to solve
\begin{equation*}\label{S6_ht_3_P1pde}
\frac{\partial P_1}{\partial\tau}=\frac{\partial^2 P_0}{\partial\xi^2}+\Big(\frac{1}{\xi}-1\Big)\,\frac{\partial P_0}{\partial\xi}+\frac{P_0}{\xi^2}-\frac{P_1}{\xi}
\end{equation*}
with the initial condition $P_1(\xi,0)=-1/\xi^2$, and hence
\begin{equation}\label{S6_ht_3_P1}
P_1(\xi,\tau)=\Big[\frac{\tau-1}{\xi^2}+\frac{4\tau-\tau^2}{2\xi^3}-\frac{3\tau^2}{2\xi^4}+\frac{\tau^3}{3\xi^5}\Big]e^{-\tau/\xi}.
\end{equation}
Using (\ref{S6_ht_3_P0}) and (\ref{S6_ht_3_P1}) in (\ref{S6_ht_3_P}) leads to (\ref{S2_th3_3}). This result is valid for $n,\,t\to\infty$ with $n/t=\xi/\tau=O(1)$.

Next, we consider the scale $n=\eta/\epsilon^2=O(\epsilon^{-2})$ and $t=\sigma/\epsilon^3=O(\epsilon^{-3})$. We note that from (\ref{S6_4_pnt}), by setting $\rho=1-\epsilon$, we have
\begin{equation*}
\rho^{-n/2}\,e^{-(1-\sqrt{\rho})^2\,t}\,=\,\exp\Big\{\frac{1}{\epsilon}\Big(\frac{\eta}{2}-\frac{\sigma}{4}\Big)+O(1)\Big\}.
\end{equation*}
This leads us to seek the asymptotic expression of $p_n(t)$ in the form
\begin{equation}\label{S6_ht_4_assume}
p_n(t)\sim \epsilon^{\nu_3}\,\widetilde{\Delta}(\eta,\sigma)\exp\Big\{\frac{1}{\epsilon}\Big[ \widetilde{\Phi}(\eta,\sigma)+\frac{\eta}{2}-\frac{\sigma}{4}\Big]\Big\}.
\end{equation}
Using (\ref{S6_ht_4_assume}) in (\ref{S2_recu1}), we have the following PDE:
\begin{equation}\label{S6_ht_4_pde}
\widetilde{\Phi}_\sigma=\widetilde{\Phi}_\eta^2-\frac{1}{\eta}.
\end{equation}
This PDE is essentially the same as that in (\ref{S6_4_pde1}). However, now we impose the initial condition $\widetilde{\Phi}(\eta,0)=-\eta/2$, which is necessary if the expansion in (\ref{S6_ht_4_assume}) is to match to (\ref{S6_ht_3_P0}) as $\tau/\xi=\sigma/(\epsilon\,\eta)\to\infty$. Thus we must solve (\ref{S6_ht_4_pde}) using characteristic curves (rays) that start from $(\eta,\sigma)=(\eta,0)$, with $\eta>0$. We again find that the geometry of the rays naturally divides the $(\eta,\sigma)$ plane into 3 parts (see Figure 2). In Figure 2 the rays in region (1) always have $d\sigma/d\eta>0$ and these never hit $\eta=0$ (the scaled time axis). Regions (2) and (3) are filled by rays that do hit $\eta=0$, and region (2) has $d\sigma/d\eta>0$ along a ray, while region (3) has $d\sigma/d\eta<0$.

The dashed curve is (\ref{S5_4_curve1}), which is a ray corresponding to $\widetilde{\Phi}_\sigma=\widetilde{A}=0$ in (\ref{S2_th3_4(1)A}), that separates regions (1) and (2). Letting $\widetilde{B}\to 1/\eta$ in (\ref{S2_th3_4(2)B}), we obtain (\ref{S5_4_curve2}), which is not a ray, and is shown by the dotted curve, which also separates regions (2) and (3). Expressions (\ref{S2_th3_4(1)phi}), (\ref{S2_th3_4(2)phi}) and (\ref{S2_th3_4(3)phi}) are obtained in the corresponding regions by solving (\ref{S6_ht_4_pde}). The function $\widetilde{\Delta}(\eta,\sigma)$ cannot be determined completely by the perturbation method, but some partial results can be obtained by using matching arguments, as was the case when $\rho<1$.

\newpage
\begin{figure}[hbp]
\begin{center}
\includegraphics[angle=0, width=0.8\textwidth]{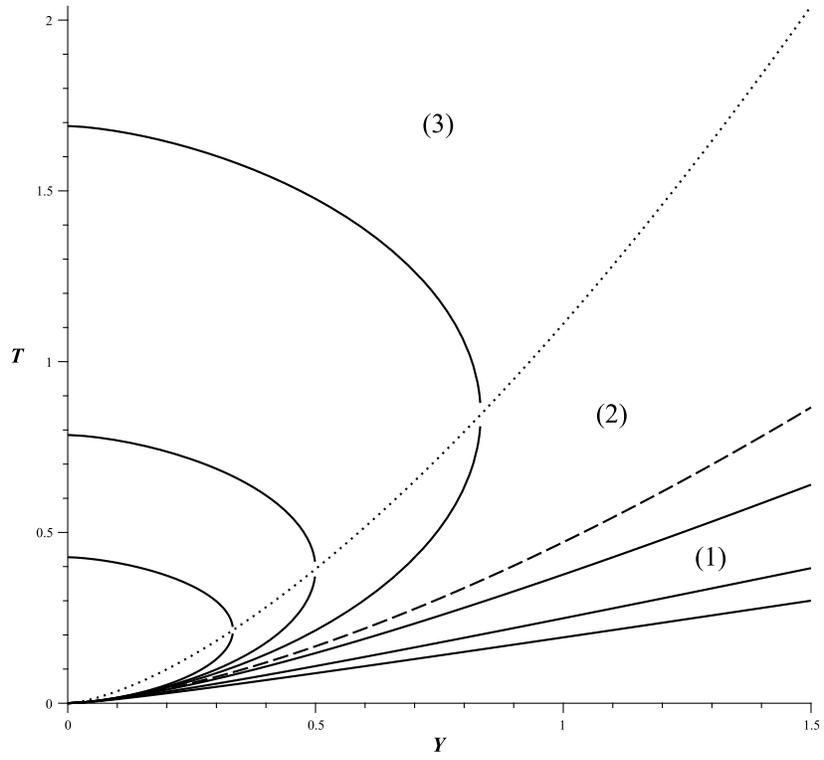}
\caption{The rays in the $(Y,T)$ plane for the $\rho<1$ case.} \label{figure1}
\end{center}
\end{figure}

\newpage
\begin{figure}[hbp]
\begin{center}
\includegraphics[angle=0, width=0.8\textwidth]{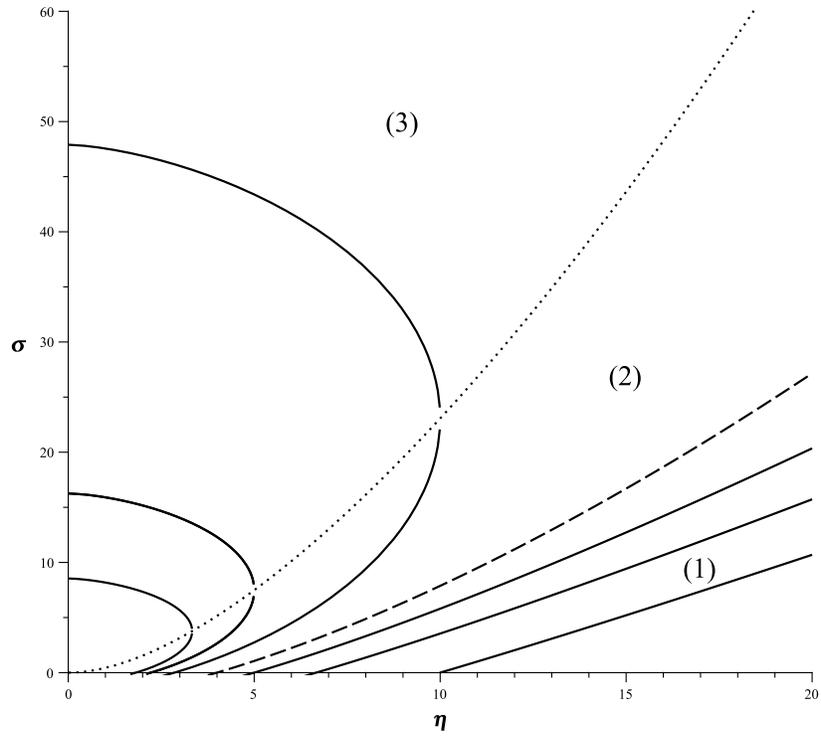}
\caption{The rays in the $(\eta,\sigma)$ plane for the heavy traffic case.} \label{figure2}
\end{center}
\end{figure}

\end{document}